\documentclass{amsart}

\pdfoutput=1

\usepackage{amsmath}
\usepackage{amssymb}
\usepackage{url}
\usepackage{mathrsfs}

\newtheorem{Thm}{Theorem}
\newtheorem{Lem}[Thm]{Lemma}
\newenvironment{Proof}{\begin{proof}}{\end{proof}}

\newcommand{\LF}{\ensuremath{\mathrm{LF}}}
\newcommand{\LFp}[3]{\ensuremath{\mathrm{LF}({#3},{#1},{#2})}}
\newcommand{\Par}{\ensuremath{\mathrm{Par}}}
\newcommand{\Is}{\ensuremath{\mathrm{Is}}}
\newcommand{\Ty}{\ensuremath{\mathrm{Ty}}}
\newcommand{\onef}[3]{\ensuremath{\mathcal{I}({#1},{#2},{#3})}}
\newcommand{\ext}{\ensuremath{\mathrm{ext}}}
\newcommand{\seeds}{\ensuremath{\mathrm{seeds}}}

\newcommand{\mydom}{\mathscr{L}}
\newcommand{\myddom}{\mathscr{L}}
\newcommand{\mydddom}{\mathscr{C}}

\newcommand{\myLatin}{\mathcal{L}}
\newcommand{\mySeed}{\mathcal{S}}

\begin{document}


\title[Latin squares of order 11]{The number of Latin squares of order 11}
\author{Alexander Hulpke}
\address{%
Department of Mathematics,
Colorado State University,
1874 Campus Delivery,
Fort Collins, CO 80523-1874, USA}
\email{hulpke@math.colostate.edu}
\author{Petteri Kaski}
\address{%
Helsinki Institute for Information Technology HIIT,
University of Helsinki, Department of Computer Science,
P.O.\ Box 68, 00014 University of Helsinki, Finland}
\email{petteri.kaski@cs.helsinki.fi}
\author{Patric R. J. \"Osterg\aa rd}
\address{%
Department of Communications and Networking,
Aalto University,
P.O.\ Box 13000, 00076 Aalto, Finland}
\email{patric.ostergard@tkk.fi}


\begin{abstract}
Constructive and nonconstructive techniques are employed 
to enumerate Latin squares and related objects.
It is established that there are
(i) $2036029552582883134196099$
main classes of Latin squares of order $11${}$;$
(ii) $6108088657705958932053657$
isomorphism classes of one-factorizations of $K_{11,11}${}$;$
(iii) $12216177315369229261482540$
isotopy classes of Latin squares of order $11${}$;$
(iv) $1478157455158044452849321016$
isomorphism classes of loops of order $11${}$;$ and 
(v) $19464657391668924966791023043937578299025$
isomorphism classes of quasigroups of order $11$.
The enumeration is constructive for the $1151666641$ 
main classes with an autoparatopy group of order at least $3$.
\end{abstract}



\maketitle


\section{Introduction}

\subsection{Latin squares}

A \emph{Latin square} of {\em order} $n$ is an $n \times n$
array over $n$ symbols such that every row and column is 
a permutation of the symbols.
Latin squares are among the most fundamental combinatorial 
objects, and have many equivalent formulations \cite[Theorem 1.11]{CDW}.

Latin squares are easily seen to exist for every $n$. 
It is, however, considerably more challenging to count the number of 
distinct Latin squares. In fact, the exact number is known only 
for $n\leq 11$, with the most recent case $n=11$ settled by McKay 
and Wanless \cite{MW}.

Here our interest is to enumerate Latin squares
up to equivalence, whereby we
focus on {\em main classes} (see \S\ref{sect:definitions}).
We refer to \cite{MMM} for a brief historical 
account of enumerating small ($n \leq 6$) Latin squares; 
this history goes as far back as the 18th century 
and includes pioneering combinatorialists such as Cayley and Euler. 
There is a unique main class for $n\leq 3$, two main classes for
$n=4,5$, and 12 main classes for $n=6$. There are 147 main
classes for $n=7$; this classification was carried out by 
Norton \cite{N} and corrected by Sade \cite{S}. 
Kolesova, Lam, and Thiel \cite{KLT} showed that there are
283657 main classes for $n=8$. All these results were constructive,
that is, representative Latin squares were (or could have been)
obtained from each main class. Recently, McKay, Meynert, and Myrvold \cite{MMM}, 
on the other hand, counted in a nonconstructive way the number of main 
classes for $n=9$ and $n=10$; these numbers are
19270853541 and 34817397894749939, respectively. 

In this paper, we enumerate the main classes for $n=11$ 
together with related objects such as quasigroups
and loops (see \cite{MMM} and \S\ref{sect:completing}). 

\begin{Thm}
\label{thm:main}
There are
\begin{enumerate}
\renewcommand{\labelenumi}{(\roman{enumi})}
\item
$2036029552582883134196099$
main classes of Latin squares of order\/ $11${}$;$
\item
$6108088657705958932053657$
isomorphism classes of one-factorizations of $K_{11,11}${}$;$
\item
$12216177315369229261482540$
isotopy classes of Latin squares of order\/ $11${}$;$
\item
$1478157455158044452849321016$
isomorphism classes of loops of order\/ $11${}$;$ and
\item
$19464657391668924966791023043937578299025$
isomorphism classes of qua\-sigroups of order\/ $11$.
\end{enumerate}
\end{Thm}

\subsection{Group actions}
\label{sect:methods}

We assume familiarity with group actions \cite{R}.
The Orbit-Stabilizer Theorem
provides the methodological foundation
for counting in \cite{MMM} as well as in this paper; 
a similar approach is also used in \cite{KO2} in the
context of one-factorizations of a complete graph.

Let $\Gamma$ be a finite group that acts on a finite set $\mydom$.
For an $X\in\mydom$, denote the orbit of $X$ by $[X]=\{\gamma X:\gamma\in\Gamma\}$
and the stabilizer of $X$ in $\Gamma$ by $\Gamma_{\!X}=\{\gamma\in\Gamma:\gamma X=X\}$.
For a subgroup $\Pi\leq\Gamma$, denote the conjugacy class of 
$\Pi$ in $\Gamma$ by $[\Pi]=\{\gamma\Pi\gamma^{-1}:\gamma\in\Gamma\}$
and the normalizer of $\Pi$ in $\Gamma$ by
$\Gamma_{\!\Pi}=\left\{\gamma\in \Gamma:\gamma\Pi\gamma^{-1}=\Pi\right\}$. 
The Orbit-Stabilizer Theorem connects the size of an orbit with
the size of a stabilizer subgroup by $|\Gamma|=|[X]|\cdot |\Gamma_{\!X}|$.

In the context of combinatorial enumeration, the elements of $\mydom$ are 
called ``labeled'' objects, and the orbits of $\Gamma$ on $\mydom$ are 
called ``unlabeled'' objects.
Our interest is to count the unlabeled objects. 
To this end, denote by $N_i$ the number of orbits 
$[X]$ with $|\Gamma_{\!X}|=i$. Because $\mydom$ partitions into orbits, 
we have
\begin{equation}
\label{eq:totaln}
|\mydom| = |\Gamma|\sum_i \frac{N_i}{i}\,.
\end{equation}
If we now know $|\Gamma|$, $|\mydom|$, and $N_i$ for 
each $i\geq 2$, then we can solve \eqref{eq:totaln} for $N_1$ 
and obtain the number of unlabeled objects $\sum_i N_i$.
In the earlier studies \cite{KO2,MMM},
the approach is to determine $|\mydom|$ by a nonconstructive 
counting technique and the values $N_i$ for $i\geq 2$ 
by a constructive classification up to isomorphism.

In this paper, we extend the nonconstructive tools to 
the setting of prescribed symmetry, which enables us to 
determine the values $N_i$ nonconstructively also for $i>1$.
There are two high-level ideas. First, 
the stabilizer subgroups of the elements of an orbit form
a conjugacy class of subgroups of $\Gamma$.
Thus, we can split the task of determining $N_i$ into subtasks 
indexed by conjugacy classes $[\Pi]$ with $|\Pi|=i$.
Second, we can solve each such subtask via an analog of \eqref{eq:totaln}. 
Namely, instead of counting individual labeled objects $X\in\mydom$ as 
in \eqref{eq:totaln}, we count formal pairs $(X,\Sigma)$ consisting 
of a labeled object $X\in\mydom$ and a subgroup $\Sigma\leq\Gamma_{\!X}$ 
with $\Sigma\in[\Pi]$.

Let us now proceed with a more detailed treatment.
Denote by $N_{[\Pi]}$ the number of orbits $[X]$ with 
$\Gamma_{\!X}\in[\Pi]$. 
Taking the sum over all conjugacy classes $[\Pi]$ with $|\Pi|=i$,
we have
\begin{equation}
\label{eq:totalni}
N_i\ =\sum_{\substack{[\Pi]:|\Pi|=i}} N_{[\Pi]}\,.
\end{equation}
For a conjugacy class $[\Pi]$, define the following two sets of formal pairs:
\begin{equation}
\label{eq:lc-def}
\begin{split}
\myddom_{[\Pi]}&=\left\{(X,\Sigma):X\in\mydom,\,\Sigma\leq\Gamma_{\!X},\,\Sigma\in[\Pi]\right\},\\[\medskipamount]
\mydddom_{[\Pi]}&=\left\{(X,\Sigma):X\in\mydom,\,\Sigma<\Gamma_{\!X},\,\Sigma\in[\Pi]\right\}.
\end{split}
\end{equation}
For $X\in\mydom$, let
\begin{equation}
\label{eq:subgroup-count}
S_{[\Pi]}(X)=\left\{\Sigma\leq\Gamma_{\!X}:\Sigma\in[\Pi]\right\}.
\end{equation}
Observe that $S_{[\Pi]}(\gamma X)=\gamma S_{[\Pi]}(X)\gamma^{-1}$
for all $\gamma\in\Gamma$.
By the Orbit--Stabilizer Theorem, we thus have
\begin{equation}
\label{eq:lc}
\left|\myddom_{[\Pi]}\right|=|\Gamma|\sum_{[X]}\,\frac{|S_{[\Pi]}(X)|}{|\Gamma_{\!X}|}\,,\qquad
\left|\mydddom_{[\Pi]}\right|=|\Gamma|\sum_{\substack{[X]:{}|\Pi|<|\Gamma_{\!X}|}}
  \frac{|S_{[\Pi]}(X)|}{|\Gamma_{\!X}|}\,.
\end{equation}
Because $\Gamma_{\!X}\in[\Pi]$ implies 
both $S_{[\Pi]}(X)=\{\Gamma_{\!X}\}$ and $|\Gamma_X|=|\Pi|$, 
we have (cf.~\eqref{eq:totaln})
\begin{equation}
\label{eq:main}
\left|\myddom_{[\Pi]}\right|=|\Gamma|\,\frac{N_{[\Pi]}}{|\Pi|}+\left|\mydddom_{[\Pi]}\right|.
\end{equation}

Equation \eqref{eq:main} is now the crux of our counting approach:
First, we compute $|\myddom_{[\Pi]}|$ nonconstructively.
Second, we constructively enumerate the orbits $[X]$ with 
$|\Pi|<|\Gamma_{\!X}|$, which gives us $|\mydddom_{[\Pi]}|$ 
via \eqref{eq:subgroup-count} and \eqref{eq:lc}.
Finally, we solve \eqref{eq:main} for $N_{[\Pi]}$.

One further simplification is that
the task of computing $|\myddom_{[\Pi]}|$ can be reduced 
from $[\Pi]$ to a representative group $\Pi\leq\Gamma$ in the 
conjugacy class. To this end, let
\begin{equation}
\myddom_{\Pi}=\left\{X\in\mydom:\Pi\leq\Gamma_{\!X}\right\}.
\end{equation}
From \eqref{eq:lc-def}, 
the fact that 
$\gamma\myddom_{\Pi}=\myddom_{\gamma\Pi\gamma^{-1}}$ 
for all $\gamma\in\Gamma$, 
and the Orbit-Stabilizer Theorem,
we have 
\begin{equation}
\label{eq:lpi-size}
|\myddom_{[\Pi]}|=\sum_{\Sigma\in[\Pi]}|\myddom_{\Sigma}|
                 =|[\Pi]|\cdot |\myddom_\Pi|
                 =\frac{|\Gamma|}{|\Gamma_{\!\Pi}|}\,|\myddom_{\Pi}|\,.
\end{equation}

We now proceed to apply this approach in the setting of Latin squares.

\subsection{Main classes}
\label{sect:definitions}

We work with the standard triple system representation
of Latin squares. To this end, let $R=\{r_1,r_2,\ldots,r_n\}$, 
$C=\{c_1,c_2,\ldots,c_n\}$, and $S=\{s_1,s_2,\ldots,s_n\}$ be three 
pairwise disjoint $n$-element sets. The intuition is that $R$ indexes 
the rows, $C$ indexes the columns, and $S$ is the set of symbols. 
We use design-theoretic terminology and call 
the elements of $R\cup C\cup S$ \emph{points} and 
the sets $R,C,S$ \emph{point classes}. 

A Latin square of order $n$ can now be represented as a set of
$n^2$ triples, where each triple $\{r_i,c_j,s_k\}$ indicates
that the symbol $s_k$ occurs in row $r_i$, column $c_j$ of the array.
Put otherwise, a {\em Latin square} (in triple system representation)
is a set $\myLatin$ of triples over $R\cup C\cup S$ such that 
(a) 
any two points from different point classes occur together 
in a unique triple; and
(b) 
no triple contains two points from the same point class. 

Assuming $R,C,S$ are fixed but arbitrary, denote by 
$\mydom$ the set of all Latin squares in triple system
representation. A ``Latin square'' in what follows 
refers to an element of $\mydom$.
Accordingly it is convenient to assume in what follows 
that a ``pair'' (of points) refers to a set of two points from
two different point classes, and a ``triple'' (of points) refers to
a set of three points from three different point classes. 

Let $\Gamma$ be the group consisting of all permutations
of $R\cup C\cup S$ that fix the partition $\{R,C,S\}$ setwise.
The structure of this group is a wreath product
$\Gamma=S_n\wr S_3\cong (S_n)^3\rtimes S_3$ with the
three copies of $S_n$ permuting the points within each of the three point classes,
and $S_3$ permuting the point classes. In particular, $|\Gamma|=3!(n!)^3$.

The orbits of $\Gamma$ on $\mydom$ are called {\em main classes}.
Two Latin squares in the same main class are {\em paratopic}.
The stabilizer of a Latin square in $\Gamma$ is called 
the {\em autoparatopy group}; the elements of the group are
{\em autoparatopisms}. 

\subsection{Proof outline}

The rest of this paper constitutes a proof of Theorem~\ref{thm:main}.
Our approach is to determine $N_i$ 
for $i=2$ via nonconstructive techniques enabled by \S\ref{sect:methods}
and for $i>2$ via constructive enumeration of main class representatives.
Finally, we solve for $N_1$ using \eqref{eq:totaln}, 
with $|\mydom|$ given by \cite{MW}. 

We begin in \S\ref{sect:autoparatopisms} by establishing the necessary 
background on autoparatopisms of Latin squares.
In \S\ref{sect:constructive}, we carry out a constructive enumeration
of all main classes of Latin squares admitting an autoparatopy group 
of order at least $3$.
In \S\ref{sect:nonconstructive}, we carry out a nonconstructive 
enumeration of Latin squares admitting fixed autoparatopisms of order $2$.
Finally, in \S\ref{sect:completing} we complete the proof of 
Theorem~\ref{thm:main} by applying 
\eqref{eq:totaln},
\eqref{eq:totalni}, 
\eqref{eq:subgroup-count},
\eqref{eq:lc},
and
\eqref{eq:main}
to the results of 
\S\ref{sect:constructive} and \S\ref{sect:nonconstructive}.

\section{Groups of autoparatopisms}
\label{sect:autoparatopisms}

\subsection{Prime order}

A group element $\alpha\in\Gamma$ of prime order $p$
is characterized up to conjugation in $\Gamma$ by 
(a) 
the prime order $p$;
(b) 
the number of points fixed in each of the point classes $R$, $C$, and $S$; and
(c) 
the number of point classes $R,C,S$ fixed.
Denote by $f_R$, $f_C$, and $f_S$ the number of points fixed 
by $\alpha$ in $R$, $C$, and $S$, respectively.
Up to conjugacy we may assume $f_R\geq f_C\geq f_S$.
Denote by $F$ the number of point classes 
$R,C,S$ fixed by $\alpha$.

We proceed to narrow down the possible types 
$(p,f_R,f_C,f_S,F)$. The results in the following two
lemmata are well known; see, for example, \cite[Theorem 1]{MMM}.

\begin{Lem}
\label{lem:autotype1}
The order $p$ divides\/ $3-F$. 
For each parameter $f$ in $f_R,f_C,f_S$ it holds
that $f=0$ or $p$ divides $n-f$; furthermore, the
latter property must hold for at least $F$
of the parameters.
\end{Lem}
\begin{Proof}
A cycle decomposition of $\alpha$ consists only of $p$-cycles 
and fixed points. So if a point class is fixed---there are
$F$ such point classes---and contains
$f$ fixed points, then $p$ divides $n-f$. Similarly,
$p$ divides $3-F$.
\end{Proof}

\begin{Lem}
\label{lem:autotype2}
Any group of autoparatopisms of a Latin square that fixes 
points from at least two point classes has the property that the
fixed points induce a subsquare.
In particular, if $f_R\geq 1$ and $f_C\geq 1$ for
a nonidentity autoparatopism, then 
$f_R=f_C=f_S\leq n/2$.
\end{Lem}
\begin{Proof}
Consider the triples in a Latin square. If a triple
has two fixed points relative to a group of autoparatopisms,
then also the third point must be fixed; otherwise the two
fixed points would occur in at least two triples, which is
impossible. Thus, the number of triples with three fixed 
points is $f_Rf_C=f_Rf_S=f_Cf_S$, from which we conclude
$f_R=f_C=f_S$ when $f_R,f_C\geq 1$. Moreover, the triples
with three fixed points induce a Latin square on the fixed points.
The order of a proper subsquare in a Latin square of order $n$ 
is at most $n/2$ \cite[Theorem 1.42]{CDW}.
\end{Proof}

We must consider the prime orders $p=2,3,5,7,11$ for $n=11$. 
Combining Lemmata \ref{lem:autotype1} and \ref{lem:autotype2}, 
we can exclude all other types $(p,f_R,f_C,f_S,F)$
except those listed in Table \ref{tbl:primeauto}.
Table~\ref{tbl:primeauto} also displays data 
related to the constructive enumeration to be discussed 
in \S\ref{sect:constructive}.

\begin{table}
\caption{Autoparatopisms of prime order}
\label{tbl:primeauto}
\begin{center}
\begin{tabular}{rrrrrrrrr}
\hline
Type &
$p$ & 
$f_R$ & 
$f_C$ & 
$f_S$ & 
$F$ &
$T$ &
Seeds &
$M_T$ \\
\hline
1 &2   & 1     & 1     & 1     & 3 &      - &       - &       - \\
2 &2   & 3     & 3     & 3     & 3 &      - &       - &       - \\
3 &2   & 5     & 5     & 5     & 3 &      - &       - &       - \\
4 &2   & 1     & 0     & 0     & 1 &      - &       - &       - \\
5 &2   & 3     & 0     & 0     & 1 &      - &       - &       - \\
6 &2   & 5     & 0     & 0     & 1 &      - &       - &       - \\
7 &2   & 7     & 0     & 0     & 1 &      - &       - &       - \\
8 &2   & 9     & 0     & 0     & 1 &      - &       - &       - \\
9 &2   & 11    & 0     & 0     & 1 &      - &       - &       - \\
10&3   & 2     & 2     & 2     & 3 &     f3 &      13 &       4 \\
11&3   & 5     & 5     & 5     & 3 &     f3 &      30 &      25 \\
12&3   & 0     & 0     & 0     & 0 &    c2n &  339638 &     165 \\
13&5   & 1     & 1     & 1     & 3 &     f3 &       7 &       1 \\
14&7   & 4     & 4     & 4     & 3 &     f3 &       6 &      16 \\
15&11  & 0     & 0     & 0     & 3 &     m1 &      28 &      33 \\
16&11  & 11    & 0     & 0     & 3 &   f1m1 &       1 &     242 \\
\hline
\end{tabular}
\end{center}
\end{table}

\subsection{Order four}

We next determine the subgroups $\Pi\leq\Gamma$ of order $4$ whose elements
of order 2 belong to the conjugacy classes
in Table \ref{tbl:primeauto}. We call such elements 
of order $2$ {\em eligible} elements.

We first compute the conjugacy
classes of elements of $\Gamma$ (there are 34048 classes) using standard
algorithms and isolate the $374$ classes of order $4$.
Amongst these we find $6$
classes for which the squares of elements are eligible.
We observe that the elements in each of these
classes are conjugate to their inverses. Therefore there are $6$ classes of 
cyclic groups of order $4$. They are listed as types 38--43 in
Table~\ref{tbl:groups4}, with generator permutations
(column ``Generators'') appearing in Table~\ref{tbl:gens4}.

We note that the conjugacy classes of eligible elements, as well as of 
elements of order $4$ with eligible squares, are determined by their 
permutation cycle structure
$1^{a_1}2^{a_2}4^{a_4}$ for $a_1+2a_2+4a_4=3n$ and 
the number $F$ of point classes $R,C,S$ that they fix. 
Table \ref{tbl:elemtype4} indicates these $16$ conjugacy classes.

\begin{table}
\caption{Autoparatopism types for order dividing 4}
\label{tbl:elemtype4}
\begin{center}
\begin{tabular}{rrr}
\hline
Type & 
Cycles & 
$F$\\
\hline
1  &$1^{3}2^{15}$    &$3$     \\
2  &$1^{9}2^{12}$    &$3$     \\
3  &$1^{15}2^{9}$    &$3$     \\
4  &$1^{1}2^{16}$    &$1$     \\
5  &$1^{3}2^{15}$    &$1$     \\
6  &$1^{5}2^{14}$    &$1$     \\
7  &$1^{7}2^{13}$    &$1$     \\
8  &$1^{9}2^{12}$    &$1$     \\
9  &$1^{11}2^{11}$   &$1$     \\
10 &$1^{3}2^{3}4^{6}$&$3$     \\
11 &$1^{5}2^{2}4^{6}$&$3$     \\
12 &$1^{7}2^{1}4^{6}$&$3$     \\
13 &$1^{9}4^{6}$     &$3$     \\
14 &$1^{1}2^{4}4^{6}$&$1$     \\
15 &$1^{3}2^{3}4^{6}$&$1$     \\
16 &$1^{33}$         &$3$     \\
\hline
\end{tabular}
\end{center}
\end{table}

To classify elementary abelian subgroups of order $4$, we note that each
such subgroup must contain three elements of order $2$, any two of which
generate the group. We will first classify such generating pairs
$(\alpha,\beta)$, where both $\alpha\in\Gamma$ and $\beta\in\Gamma$ are 
eligible, and $\beta$ is contained in the centralizer of $\alpha$, which we
denote by
$\Gamma_\alpha=\{\gamma\in\Gamma:\gamma\alpha\gamma^{-1}=\alpha\}$. 
Furthermore
$\alpha\beta$ must be eligible as well. Clearly, any two of the
three elements $\alpha$, $\beta$, $\alpha\beta$ generate the same
group. We therefore can assume without loss of generality that 
(in an arbitrary total order on the conjugacy classes of $\Gamma$, for
example, the ``Type'' number in Table \ref{tbl:elemtype4}) the class of 
$\alpha$ is not larger than the class of $\beta$, which is not larger than 
the class of $\alpha\beta$. If this is fulfilled, we call the pair
$(\alpha,\beta)$ {\em valid}.

Furthermore, we need to classify such pairs only up to
conjugacy in $\Gamma$. For this, we may assume that $\alpha$ is the chosen
representative of its conjugacy class and $\beta$ is chosen up to conjugacy by
$\Gamma_\alpha$. We determine the representatives for $\beta$ by 
computing for each eligible $\alpha$ the conjugacy classes of 
$\Gamma_\alpha$ and identifying
amongst these the classes of eligible elements. For each such class
representative $\beta$ we test whether $(\alpha,\beta)$ is valid, and 
collect the valid pairs in a list.

Because conjugacy of subgroups conjugates elements to elements, the
lexicographic choice of
generating pairs ensures that the only remaining potential
conjugacy amongst subgroups
generated by such pairs can happen if both $\alpha$ and $\beta$, 
or $\beta$ and $\alpha\beta$ are in the same class under $\Gamma$.

A classification of such pairs on the computer finds $37$ valid pairs; none
of the corresponding subgroups are conjugate, as verified by
explicit conjugacy tests in $\Gamma$.
These groups are listed as types 1--37 in Table~\ref{tbl:groups4}, again
generators appear in Table~\ref{tbl:gens4}.

Calculations were done in the system {\sf GAP} \cite{GAP}.

Up to conjugacy in $\Gamma$, there are therefore exactly 
43 groups $\Pi\leq\Gamma$ of order $4$ whose 
elements of order 2 belong to the conjugacy classes
in Table \ref{tbl:primeauto}. 

It turns out that
these groups are distinguished up to conjugacy in $\Gamma$ by the conjugacy 
classes of their elements (column ``Elements'')
and the lengths $1^{a_1}2^{a_2}4^{a_4}$ for $a_1+2a_2+4a_4=3n$
of their orbits on $R\cup C\cup S$ (column ``Orbits'').

\begin{table}
\caption{Candidate groups of order 4}
\label{tbl:groups4}
\begin{center}
{\small
\begin{tabular}{rrrrrrr}
\hline
Type & 
Generators &
Elements & 
Orbits &
$T$ &
Seeds &
$M_T$ 
\\
\hline
1*\!\!\! & $\gamma_{1},\gamma_{2}$ &3, 6, 9, 16&    $1^{5}2^{8}4^{3}$   & -&-&-\\
2 & $\gamma_{1},\gamma_{3}$ &3, 4, 7, 16&    $1^{1}2^{10}4^{3}$  & f1m2 &     102 &  11 \\
3*\!\!\! & $\gamma_{1},\gamma_{4}$ &2, 3, 3, 16&    $1^{3}2^{15}$       & -&-&-\\
4*\!\!\! & $\gamma_{1},\gamma_{5}$ &2, 3, 3, 16&    $1^{7}2^{9}4^{2}$   & -&-&-\\
5*\!\!\! & $\gamma_{1},\gamma_{6}$ &3, 7, 8, 16&    $1^{5}2^{8}4^{3}$   & -&-&-\\
6 & $\gamma_{1},\gamma_{7}$ &3, 5, 6, 16&    $1^{1}2^{10}4^{3}$  & f1m2 &     230 &  11 \\ 
7 & $\gamma_{1},\gamma_{8}$ &1, 2, 3, 16&    $1^{3}2^{9}4^{3}$   & f3   &      40 &   1 \\ 
8 & $\gamma_{1},\gamma_{9}$ &3, 6, 7, 16&    $1^{5}2^{6}4^{4}$   & f3n  &      80 &  10 \\ 
9*\!\!\! & $\gamma_{1},\gamma_{10}$&3, 4, 5, 16&    $1^{1}2^{8}4^{4}$   & -&-&-\\
10*\!\!\!& $\gamma_{1},\gamma_{11}$&2, 3, 3, 16&    $1^{5}2^{12}4^{1}$  & -&-&-\\
11& $\gamma_{1},\gamma_{12}$&3, 5, 8, 16&    $1^{3}2^{9}4^{3}$   & f3n  &      80 &   1 \\ 
12& $\gamma_{1},\gamma_{13}$&3, 6, 7, 16&    $1^{3}2^{9}4^{3}$   & f3n  &      80 &   1 \\ 
13*\!\!\!& $\gamma_{1},\gamma_{14}$&2, 3, 3, 16&    $1^{9}2^{6}4^{3}$   & -&-&-\\
14& $\gamma_{1},\gamma_{15}$&3, 5, 6, 16&    $1^{3}2^{7}4^{4}$   & f3n  &      80 &   1 \\ 
15&$\gamma_{16},\gamma_{17}$&2, 2, 2, 16&    $1^{9}4^{6}$        & f3   &      12 &   9 \\ 
16*\!\!\!&$\gamma_{16},\gamma_{18}$&2, 2, 2, 16&    $1^{5}2^{6}4^{4}$   & -&-&-\\
17&$\gamma_{16},\gamma_{19}$&1, 1, 2, 16&    $1^{3}2^{3}4^{6}$   & f3   &      20 &   1 \\ 
18&$\gamma_{16},\gamma_{20}$&2, 5, 9, 16&    $1^{3}2^{7}4^{4}$   & f3n  &     127 &   1 \\ 
19&$\gamma_{16},\gamma_{21}$&2, 7, 7, 16&    $1^{3}2^{7}4^{4}$   & f3n  &      83 &   1 \\ 
20&$\gamma_{16},\gamma_{22}$&2, 5, 5, 16&    $1^{3}2^{3}4^{6}$   & f3n  &      83 &   1 \\ 
21&$\gamma_{16},\gamma_{23}$&2, 6, 8, 16&    $1^{3}2^{7}4^{4}$   & f3n  &     127 &   1 \\ 
22&$\gamma_{16},\gamma_{24}$&2, 4, 8, 16&    $1^{1}2^{8}4^{4}$   & f1m2 &     608 &  11 \\ 
23&$\gamma_{16},\gamma_{25}$&2, 6, 6, 16&    $1^{1}2^{8}4^{4}$   & f1m1 &    2215 &  22 \\ 
24*\!\!\!&$\gamma_{16},\gamma_{26}$&2, 4, 4, 16&    $1^{1}2^{4}4^{6}$   & -&-&-\\
25&$\gamma_{16},\gamma_{27}$&2, 5, 7, 16&    $1^{1}2^{8}4^{4}$   & f1m1 &    2835 &  22 \\ 
26&$\gamma_{16},\gamma_{28}$&2, 2, 2, 16&    $1^{3}2^{9}4^{3}$   & f3   &      10 &   1 \\ 
27*\!\!\!&$\gamma_{16},\gamma_{29}$&2, 2, 2, 16&    $1^{7}2^{3}4^{5}$   & -&-&-\\
28&$\gamma_{16},\gamma_{30}$&2, 5, 7, 16&    $1^{3}2^{5}4^{5}$   & f3n  &     127 &   1 \\ 
29&$\gamma_{16},\gamma_{31}$&2, 6, 6, 16&    $1^{3}2^{5}4^{5}$   & f3n  &      83 &   1 \\ 
30&$\gamma_{16},\gamma_{32}$&2, 4, 6, 16&    $1^{1}2^{6}4^{5}$   & f1m1 &    2854 &  22 \\ 
31&$\gamma_{16},\gamma_{33}$&2, 5, 5, 16&    $1^{1}2^{6}4^{5}$   & f1m1 &    3351 &  22 \\ 
32&$\gamma_{34},\gamma_{35}$&1, 4, 9, 16&    $1^{1}2^{6}4^{5}$   & f1m1 &     194 &  22 \\ 
33&$\gamma_{34},\gamma_{36}$&1, 6, 7, 16&    $1^{1}2^{6}4^{5}$   & f1m1 &    2300 &  22 \\ 
34&$\gamma_{34},\gamma_{37}$&1, 5, 8, 16&    $1^{1}2^{6}4^{5}$   & f1m1 &    1083 &  22 \\ 
35&$\gamma_{34},\gamma_{38}$&1, 4, 7, 16&    $1^{1}2^{4}4^{6}$   & f1m1 &    1696 &  22 \\ 
36&$\gamma_{34},\gamma_{39}$&1, 5, 6, 16&    $1^{1}2^{4}4^{6}$   & f1m1 &    6406 &  22 \\ 
37&$\gamma_{34},\gamma_{40}$&1, 4, 5, 16&    $1^{1}2^{2}4^{7}$   & f1m1 &    1761 &  22 \\ 
38&$\gamma_{41}$            &2, 13, 13, 16&    $1^{9}4^{6}$      & f3   &      20 &   9 \\ 
39&$\gamma_{42}$            &2, 10, 10, 16&    $1^{3}2^{3}4^{6}$ & f3   &      20 &   1 \\ 
40*\!\!\!&$\gamma_{43}$            &2, 11, 11, 16&    $1^{5}2^{2}4^{6}$ & -&-&-\\
41*\!\!\!&$\gamma_{44}$            &2, 12, 12, 16&    $1^{7}2^{1}4^{6}$ & -&-&-\\
42&$\gamma_{45}$            &2, 15, 15, 16&    $1^{3}2^{3}4^{6}$ & f3n  &       4 &   1 \\ 
43&$\gamma_{46}$            &2, 14, 14, 16&    $1^{1}2^{4}4^{6}$ & f1m1 &     242 &  22 \\ 
\hline\\[-2mm]
\makebox[0pt][l]{\hspace*{-8mm}*) Excluded by Lemma~\ref{lem:group4}}
\end{tabular}
}
\end{center}
\end{table}

\begin{table}
\caption{Generators for Table~\ref{tbl:groups4}}
\label{tbl:gens4}
\vspace*{-5mm}
{\tiny
\begin{align*}
\rule{5mm}{0.5pt}&\rule{\textwidth}{0.5pt}\\
R &= \{1,2,3,4,5,6,7,8,9,10,11\}\\
C &= \{12,13,14,15,16,17,18,19,20,21,22\}\\
S &= \{23,24,25,26,27,28,29,30,31,32,33\}\\
\rule{5mm}{0.5pt}&\rule{\textwidth}{0.5pt}\\
\gamma_{1} &= ( 1,\! 2)( 3,\! 4)( 5,\! 6)(12,\!13)(14,\!15)(16,\!17)(23,\!24)(25,\!26)(27,\!28)
\\
\gamma_{2} &= (12,\!23)(13,\!24)(14,\!25)(15,\!26)(16,\!27)(17,\!28)(18,\!31)(19,\!29)(20,\!30)
(21,\!32)(22,\!33)
\\
\gamma_{3} &= ( 8,\!10)( 9,\!11)(12,\!23)(13,\!24)(14,\!25)(15,\!26)(16,\!27)(17,\!28)(18,\!31)
(19,\!29)(20,\!30)(21,\!32)(22,\!33)
\\
\gamma_{4} &= ( 1,\! 2)( 8,\!10)( 9,\!11)(16,\!17)(19,\!22)(20,\!21)(27,\!28)(29,\!33)(30,\!32)
\\
\gamma_{5} &= ( 1,\! 2)( 8,\!10)( 9,\!11)(12,\!14)(13,\!15)(19,\!22)(23,\!25)(24,\!26)(29,\!33)
\\
\gamma_{6} &= ( 1,\! 2)(12,\!23)(13,\!24)(14,\!25)(15,\!26)(16,\!27)(17,\!28)(18,\!31)(19,\!29)
(20,\!30)(21,\!32)(22,\!33)
\\
\gamma_{7} &= ( 1,\! 2)( 8,\!10)( 9,\!11)(12,\!23)(13,\!24)(14,\!25)(15,\!26)(16,\!27)(17,\!28)
(18,\!31)(19,\!29)(20,\!30)(21,\!32)(22,\!33)
\\
\gamma_{8} &= ( 3,\! 5)( 4,\! 6)( 8,\!10)( 9,\!11)(12,\!14)(13,\!15)(19,\!22)(20,\!21)(23,\!25)
(24,\!26)(29,\!33)(30,\!32)
\\
\gamma_{9} &= ( 3,\! 5)( 4,\! 6)(12,\!23)(13,\!24)(14,\!25)(15,\!26)(16,\!27)(17,\!28)(18,\!31)
(19,\!29)(20,\!30)(21,\!32)(22,\!33)
\\
\gamma_{10} &= ( 1,\! 2)( 3,\! 5)( 4,\! 6)( 8,\!10)( 9,\!11)(12,\!23)(13,\!24)(14,\!25)(15,\!26)
(16,\!27)(17,\!28)(18,\!31)(19,\!29)(20,\!30)(21,\!32)(22,\!33)
\\
\gamma_{11} &= ( 1,\! 2)( 8,\!10)( 9,\!11)(16,\!17)(19,\!22)(20,\!21)(23,\!25)(24,\!26)(29,\!33)
\\
\gamma_{12} &= ( 9,\!11)(12,\!23)(13,\!24)(14,\!25)(15,\!26)(16,\!27)(17,\!28)(18,\!31)(19,\!29)
(20,\!30)(21,\!32)(22,\!33)
\\
\gamma_{13} &= ( 1,\! 2)( 9,\!11)(12,\!23)(13,\!24)(14,\!25)(15,\!26)(16,\!27)(17,\!28)(18,\!31)
(19,\!29)(20,\!30)(21,\!32)(22,\!33)
\\
\gamma_{14} &= ( 3,\! 5)( 4,\! 6)( 9,\!11)(12,\!14)(13,\!15)(19,\!22)(23,\!25)(24,\!26)(29,\!33)
\\
\gamma_{15} &= ( 3,\! 5)( 4,\! 6)( 9,\!11)(12,\!23)(13,\!24)(14,\!25)(15,\!26)(16,\!27)(17,\!28)
(18,\!31)(19,\!29)(20,\!30)(21,\!32)(22,\!33)
\\
\gamma_{16} &= ( 1,\! 2)( 3,\! 4)( 5,\! 6)( 7,\! 8)(12,\!13)(14,\!15)(16,\!17)(18,\!19)(23,\!24)
(25,\!26)(27,\!28)(29,\!30)
\\
\gamma_{17} &= ( 1,\! 3)( 2,\! 4)( 5,\! 7)( 6,\! 8)(12,\!14)(13,\!15)(16,\!18)(17,\!19)(23,\!25)
(24,\!26)(27,\!29)(28,\!30)
\\
\gamma_{18} &= ( 1,\! 3)( 2,\! 4)( 5,\! 7)( 6,\! 8)(12,\!15)(13,\!14)(16,\!17)(21,\!22)(23,\!26)
(24,\!25)(27,\!28)(32,\!33)
\\
\gamma_{19} &= ( 1,\! 3)( 2,\! 4)( 5,\! 7)( 6,\! 8)(10,\!11)(12,\!14)(13,\!15)(16,\!18)(17,\!19)
(21,\!22)(23,\!25)(24,\!26)(27,\!29)(28,\!30)(32,\!33)
\\
\gamma_{20} &= (12,\!23)(13,\!24)(14,\!25)(15,\!26)(16,\!27)(17,\!28)(18,\!29)(19,\!30)(20,\!31)
(21,\!32)(22,\!33)
\\
\gamma_{21} &= ( 5,\! 6)( 7,\! 8)(12,\!23)(13,\!24)(14,\!25)(15,\!26)(16,\!27)(17,\!28)(18,\!29)
(19,\!30)(20,\!31)(21,\!32)(22,\!33)
\\
\gamma_{22} &= ( 1,\! 3)( 2,\! 4)( 5,\! 7)( 6,\! 8)(12,\!23)(13,\!24)(14,\!25)(15,\!26)(16,\!27)
(17,\!28)(18,\!29)(19,\!30)(20,\!31)(21,\!32)(22,\!33)
\\
\gamma_{23} &= ( 7,\! 8)(12,\!24)(13,\!23)(14,\!26)(15,\!25)(16,\!28)(17,\!27)(18,\!29)(19,\!30)
(20,\!31)(21,\!32)(22,\!33)
\\
\gamma_{24} &= (10,\!11)(12,\!23)(13,\!24)(14,\!25)(15,\!26)(16,\!27)(17,\!28)(18,\!29)(19,\!30)
(20,\!31)(21,\!33)(22,\!32)
\\
\gamma_{25} &= ( 5,\! 6)( 7,\! 8)(10,\!11)(12,\!23)(13,\!24)(14,\!25)(15,\!26)(16,\!27)(17,\!28)
(18,\!29)(19,\!30)(20,\!31)(21,\!33)(22,\!32)
\\
\gamma_{26} &= ( 1,\! 3)( 2,\! 4)( 5,\! 7)( 6,\! 8)(10,\!11)(12,\!23)(13,\!24)(14,\!25)(15,\!26)
(16,\!27)(17,\!28)(18,\!29)(19,\!30)(20,\!31)(21,\!33)(22,\!32)
\\
\gamma_{27} &= ( 7,\! 8)(10,\!11)(12,\!24)(13,\!23)(14,\!26)(15,\!25)(16,\!28)(17,\!27)(18,\!29)
(19,\!30)(20,\!31)(21,\!33)(22,\!32)
\\
\gamma_{28} &= ( 1,\! 3)( 2,\! 4)( 7,\! 8)(10,\!11)(12,\!15)(13,\!14)(16,\!17)(21,\!22)(23,\!25)
(24,\!26)(29,\!30)(32,\!33)
\\
\gamma_{29} &= ( 1,\! 3)( 2,\! 4)( 7,\! 8)(10,\!11)(12,\!14)(13,\!15)(16,\!18)(17,\!19)(23,\!26)
(24,\!25)(27,\!30)(28,\!29)
\\
\gamma_{30} &= ( 1,\! 3)( 2,\! 4)(12,\!25)(13,\!26)(14,\!23)(15,\!24)(16,\!27)(17,\!28)(18,\!29)
(19,\!30)(20,\!31)(21,\!32)(22,\!33)
\\
\gamma_{31} &= ( 1,\! 3)( 2,\! 4)( 7,\! 8)(12,\!26)(13,\!25)(14,\!24)(15,\!23)(16,\!28)(17,\!27)
(18,\!29)(19,\!30)(20,\!31)(21,\!32)(22,\!33)
\\
\gamma_{32} &= ( 1,\! 3)( 2,\! 4)(10,\!11)(12,\!25)(13,\!26)(14,\!23)(15,\!24)(16,\!27)(17,\!28)
(18,\!29)(19,\!30)(20,\!31)(21,\!33)(22,\!32)
\\
\gamma_{33} &= ( 1,\! 3)( 2,\! 4)( 7,\! 8)(10,\!11)(12,\!26)(13,\!25)(14,\!24)(15,\!23)(16,\!28)
(17,\!27)(18,\!29)(19,\!30)(20,\!31)(21,\!33)(22,\!32)
\\
\gamma_{34} &= ( 1,\! 2)( 3,\! 4)( 5,\! 6)( 7,\! 8)( 9,\!10)(12,\!13)(14,\!15)(16,\!17)(18,\!19)
(20,\!21)(23,\!24)(25,\!26)(27,\!28)(29,\!30)(31,\!32)
\\
\gamma_{35} &= (12,\!23)(13,\!24)(14,\!25)(15,\!26)(16,\!29)(17,\!30)(18,\!27)(19,\!28)(20,\!31)
(21,\!32)(22,\!33)
\\
\gamma_{36} &= ( 7,\! 8)( 9,\!10)(12,\!23)(13,\!24)(14,\!25)(15,\!26)(16,\!29)(17,\!30)(18,\!27)
(19,\!28)(20,\!31)(21,\!32)(22,\!33)
\\
\gamma_{37} &= ( 3,\! 4)(12,\!23)(13,\!24)(14,\!25)(15,\!26)(16,\!29)(17,\!30)(18,\!27)(19,\!28)
(20,\!31)(21,\!32)(22,\!33)
\\
\gamma_{38} &= ( 7,\! 9)( 8,\!10)(12,\!23)(13,\!24)(14,\!25)(15,\!26)(16,\!29)(17,\!30)(18,\!27)
(19,\!28)(20,\!31)(21,\!32)(22,\!33)
\\
\gamma_{39} &= ( 3,\! 4)( 7,\! 9)( 8,\!10)(12,\!23)(13,\!24)(14,\!25)(15,\!26)(16,\!29)(17,\!30)
(18,\!27)(19,\!28)(20,\!31)(21,\!32)(22,\!33)
\\
\gamma_{40} &= ( 1,\! 2)( 3,\! 5)( 4,\! 6)( 7,\! 9)( 8,\!10)(12,\!23)(13,\!24)(14,\!25)(15,\!26)
(16,\!29)(17,\!30)(18,\!27)(19,\!28)(20,\!31)(21,\!32)(22,\!33)
\\
\gamma_{41} &= ( 1,\! 2,\! 3,\! 4)( 5,\! 6,\! 7,\! 8)(12,\!13,\!14,\!15)(16,\!17,\!18,\!19)
(23,\!24,\!25,\!26)(27,\!28,\!29,\!30)
\\
\gamma_{42} &= ( 1,\! 2,\! 3,\! 4)( 5,\! 6,\! 7,\! 8)( 9,\!10)(12,\!13,\!14,\!15)(16,\!17,\!18,\!19)
(20,\!21)(23,\!24,\!25,\!26)(27,\!28,\!29,\!30)(31,\!32)
\\
\gamma_{43} &= ( 1,\! 2,\! 3,\! 4)( 5,\! 6,\! 7,\! 8)(12,\!13,\!14,\!15)(16,\!17,\!18,\!19)(20,\!21)
(23,\!24,\!25,\!26)(27,\!28,\!29,\!30)(31,\!32)
\\
\gamma_{44} &= ( 1,\! 2,\! 3,\! 4)( 5,\! 6,\! 7,\! 8)(12,\!13,\!14,\!15)(16,\!17,\!18,\!19)
(23,\!24,\!25,\!26)(27,\!28,\!29,\!30)(31,\!32)
\\
\gamma_{45} &= ( 1,\!13,\! 2,\!12)( 3,\!15,\! 4,\!14)( 5,\!17,\! 6,\!16)( 7,\!19,\! 8,\!18)( 9,\!20)
(10,\!21)(11,\!22)(23,\!24,\!25,\!26)(27,\!28,\!29,\!30)
\\
\gamma_{46} &= ( 1,\!13,\! 2,\!12)( 3,\!15,\! 4,\!14)( 5,\!17,\! 6,\!16)( 7,\!19,\! 8,\!18)( 9,\!20)
(10,\!21)(11,\!22)(23,\!24,\!25,\!26)(27,\!28,\!29,\!30)(31,\!32)
\\
\rule{5mm}{0.5pt}&\rule{\textwidth}{0.5pt}
\end{align*}
}
\end{table}

We can rule out the following groups with combinatorial arguments.

\begin{Lem}
\label{lem:group4}
The groups of types\/ 
$1$, $3$, $4$, $5$, $9$, $10$, $13$, $16$, $24$, $27$, $40$, and\/ $41$ 
are not admitted by a Latin square of order\/ $11$.
\end{Lem}

\begin{proof}
Let $R = \{1,2,\ldots,11\}$,
$C = \{12,13,\ldots,22\}$, and
$S = \{23,24,\ldots,33\}$.

{\em Types\/ $4$, $10$, $16$, $27$, $40$, and\/ $41$}:
All point classes contain fixed points, but their number varies
between point classes. This contradicts Lemma~\ref{lem:autotype2}.

{\em Type\/ $1$}:
Let
\begin{alignat*}{2}
R_1 &= \{1,2,3,4,5,6\}\,,\quad &R_2 &= \{7,8,9,10,11\}\,,\\
C_1 &= \{12,13,14,15,16,17\}\,,\quad &C_2 &= \{18,19,20,21,22\}\,,\\
S_1 &= \{23,24,25,26,27,28\}\,,\quad &S_2 &= \{29,30,31,32,33\}\,.
\end{alignat*}
Lemma \ref{lem:autotype2} applied to $\gamma_1$ reveals that 
$R_2 \cup C_2 \cup S_2$ induces a $5 \times 5$ subsquare.
Thus, there exists no triple with exactly two points from 
$R_2 \cup C_2 \cup S_2$.

Next consider pairs with one point from $C_1$ 
and one from $S_1$; there are $|C_1|\cdot |S_1|=36$ such pairs.
The action of the group partitions these pairs combined with
any points from $R$ into orbits,
exactly six of which have length 4; 
for example, one such orbit is  
\[
\{\{x,12,25\},\{x,14,23\},\{\gamma_1(x),13,26\},\{\gamma_1(x),15,24\}\}\,.
\]
As $\gamma_1(x)=x$ when $x \in R_2$, 
any such $x$ can be in at most one orbit (otherwise
it would occur in at least $2 \cdot 4 = 8$ triples 
together with points from $C_1$ but $|C_1|=6<8$).

It follows that $x \in R_2$ in at most five 
of the above mentioned six orbits of length 4, whereby
$x \in R_1$ in at least one of the orbits.
Such an orbit with $x \in R_1$ contains
two pairs $\{x,y\}$ with $y \in C_1$ (and recall that the
third point of the triple is from $S_1$). Finally, consider
how pairs $\{x,z\}$ with $z \in S_2$ are covered. The third
point in the corresponding triple has to come from
$C_1$ (by the initial comment). However, we have 5
points in $S_2$ to be combined with $x$ and $6-2=4$ points
in $C_1$, which is not possible.

{\em Type\/ $3$}:
It follows by Lemma \ref{lem:autotype2} that there are 
$5 \times 5$ subsquares over
$\{7,8,9,10,11\}$, $\{18,19,20,21,22\}$, $\{29,30,31,32,33\}$ and
$\{3,4,5,6,7\}$, $\{12,13,14,15,18\}$, $\{23,24,\linebreak 25,26,31\}$
for $\gamma_1$ and $\gamma_4$, respectively.

Consider the four triples containing the point $3$ and one point
from $\{19,20,\linebreak 21,22\}$. By the given subsquare structure,
the only possibilities for the third point of the triple
are the points $\{27,28\}$, implying that $\{3,27\}$ or
$\{3,28\}$ would occur in more than one triple.

{\em Type\/ $5$}:
The proof for Type 1 is applicable as such, since the only 
difference between these cases is the additional transposition
$(1,2)$ in generator $\gamma_6$, which leads to different orbits---for
example, 
$\{1,12,25\}, \{2,14,23\}, \{2,13,26\}, \{1,15,24\}$ instead of
$\{1,12,25\}, \{1,14,23\}, \{2,13,26\}, \{2,15,24\}$---but has
no impact on the arguments of the proof.

{\em Type\/ $9$}:
The only point fixed by $\gamma_{10}$ is $7$, which is then the only
possible third point in triples containing the points
in transpositions of $\gamma_{10}$ with 
one point each from $C$ and $S$. Since there are $11$
such triples, the point $7$ does not occur in any other triples.

By considering possible choices for $x$
in $\{x,12,24\}$, we get up to symmetry four
cases with the orbits
\begin{align*}
&\{\{1,12,24\},\{2,13,23\}\}\,,\\
&\{\{3,12,24\},\{4,13,23\},\{5,13,23\},\{6,12,24\}\}\,,\\
&\{\{8,12,24\},\{8,13,23\},\{10,13,23\},\{10,12,24\}\}\,,\text{ and}\\
&\{\{7,12,24\},\{7,13,23\}\}\,.
\end{align*}
Because of the earlier
argument, we cannot have $x=7$, and because pairs of
points occur twice in two other cases, we must have
$x \in \{1,2\}$.

Repeated use of the previous reasoning gives that 
each of the pairs $\{12,24\}$, $\{13,23\}$, 
$\{14,26\}$, $\{15,25\}$, $\{16,28\}$, and $\{17,27\}$
must occur in triples together with the points $1$ and $2$
(three times each). Now consider triples 
$\{1,x,y\}$ with $y \in \{29,30,31,32,33\} = S_2$. It follows
from the subsquare implied by $\gamma_1$ that
we cannot have $x \in \{18,19,20,21,22\} = C_2$. There are
therefore $|C \setminus C_2| = 6$ possible choices for $x$, 
3 of which already occur in triples together with the point $1$. 
But $6-3 = 3 < 5 = |S_2|$.

{\em Type\/ $13$}:
By applying Lemma \ref{lem:autotype2} to $\gamma_1$ and $\gamma_1\gamma_{14}$ 
we get that there is a $5\times 5$ subsquare over $\{7,8,9,10,11\}$, 
$\{18,19,20,21,22\}$, $\{29,30,31,32,33\}$ and a 
$3\times 3$ subsquare over $\{7,8,10\}$, $\{18,20,21\}$,
$\{30,31,32\}$, respectively. But the latter is a subsquare of the
former, and a $5 \times 5$ square cannot contain a $3\times 3$ subsquare \cite[Theorem 1.42]{CDW}.

{\em Type\/ $24$}:
Analogously to the case of Type 9, we argue that,
as the only point fixed by $\gamma_{26}$ is $9$, 
this point occurs exactly in triples with the other two points 
from transpositions of $\gamma_{26}$. 

By considering possible choices for $x$
in $\{x,12,24\}$, we get up to symmetry three
cases with the orbits
\begin{align*}
&\{\{1,12,24\},\{3,13,23\},\{4,12,24\},\{2,13,23\}\}\,,\\
&\{\{10,12,24\},\{10,13,23\},\{11,13,23\},\{11,12,24\}\}\,,\text{ and}\\
&\{\{9,12,24\},\{9,13,23\}\}\,.
\end{align*}
But because of the earlier argument, we cannot have $x=9$, 
and pairs of points occur twice in the other two cases.
\end{proof}

We remark that the combinatorial proofs in 
Lemma~\ref{lem:group4} were triggered by computer
searches showing nonexistence.

\section{Constructive enumeration}
\label{sect:constructive}

\subsection{Admissible groups}

We recall the following corollary of the Sylow Theorems.

\begin{Thm}
\label{thm:ppgroup}
{\rm\bf\cite[Corollary~4.15]{R}}
Let $\Lambda$ be a finite group, let $p$ be a prime, 
and let $k$ be a positive integer.
If $p^k$ divides $|\Lambda|$, 
then $\Lambda$ contains a subgroup of order $p^k$.
\end{Thm}

Let us say that a Latin square $\myLatin$ {\em admits} 
$\Lambda\leq\Gamma$ as a group of autoparatopisms 
if $\Lambda\leq\Gamma_{\!\myLatin}$.
Combining Theorem~\ref{thm:ppgroup} with the analysis in 
\S\ref{sect:autoparatopisms}, we have that every
Latin square $\myLatin$ of order $11$ 
with $|\Gamma_{\!\myLatin}|\geq 3$ admits at least 
one subgroup conjugate to
(a) a group of odd prime order in Table~\ref{tbl:primeauto}; 
or
(b) a group of order 4 in Table~\ref{tbl:groups4}
not excluded by Lemma~\ref{lem:group4}.
Let us call these groups the {\em admissible groups}.

\subsection{The enumeration}

We apply the framework in \cite{K} to construct exactly one
representative from each main class of Latin squares
admitting at least one admissible group.
 
Associated with each conjugacy class $[\Pi]$ of admissible groups
is a set of {\em seeds}, at least one of which is guaranteed
to occur in every Latin square admitting $\Pi$.

The formal definition of a seed is as follows.
Let $\Pi\leq\Gamma$ be an admissible group. Let $T\subseteq R\cup C\cup S$
be a set of points whose size and composition is determined 
by column ``$T$'' in Tables \ref{tbl:primeauto} and \ref{tbl:groups4}
as follows:
\begin{itemize}
\item[f$i$\ :]
  indicates $i$ points fixed by $\Pi$;
\item[m$i$\ :] 
  indicates $i$ points moved by $\Pi$;
\item[c$i$\ :] 
  indicates $i$ points in same point class;
\item[n\ :]
  see Item (d) below.
\end{itemize}
(For example, ``f1m2'' indicates that $T$ consists of 3 points,
1 point fixed by $\Pi$, and 2 points moved by $\Pi$.)
Finally, let $\mySeed$ be a union of $\Pi$-orbits of triples
such that, referring to the elements of $\mySeed$ as {\em blocks}, 
\begin{itemize}
\item[(a)]
any pair occurs in at most one block;
\item[(b)]
each point in $T$ occurs in exactly $n$ blocks;
\item[(c)]
$T$ has nonempty intersection with at least one block on
every $\Pi$-orbit on $\mySeed$; and
\item[(d)]
the set $T$ occurs in at least one block {\em unless} 
the composition of $T$ has the ``n'' indicator.
\end{itemize}
Each such tuple $(\Pi,T,\mySeed)$ is a {\em seed} associated
with the conjugacy class $[\Pi]$. 

Let $\gamma\in\Gamma$ act on a seed by 
$\gamma(\Pi,T,\mySeed)=(\gamma\Pi\gamma^{-1},\gamma T,\gamma \mySeed)$.
The orbits of $\Gamma$ on the set of all seeds are 
the {\em isomorphism classes} of seeds. A Latin square $\myLatin$ 
{\em contains} (or {\em extends}) a seed $(\Pi,T,\mySeed)$ 
if $\Pi\leq\Gamma_{\!\myLatin}$ and $\mySeed\subseteq\myLatin$.

Given a seed $(\Pi,T,\mySeed)$, the task of finding all Latin squares
that extend the seed is an instance of the {\em exact cover problem}.
In particular, the task is to cover exactly once all uncovered pairs 
using $\Pi$-orbits of triples; each triple covers the pairs occurring 
in it.

Our constructive enumeration approach now proceeds as follows.
First, we classify the seeds up to isomorphism using the algorithms
described in \cite{K}. The number of nonisomorphic seeds associated 
with each conjugacy class of admissible groups is given in column ``Seeds''
in Tables \ref{tbl:primeauto} and \ref{tbl:groups4}. Once the seeds
have been classified, we use {\em libexact}\/ \cite{KP} to find all 
extensions of each seed to a Latin square. As each extension  
$\myLatin$ of a seed $(\Pi,T,\mySeed)$ is constructed, it is subjected 
to isomorph rejection tests developed in \cite{K}, which derive
from the canonical augmentation technique of McKay \cite{M3}. 
First, we identify a canonical $\Gamma_{\!\myLatin}$-orbit of seeds 
contained in $\myLatin$. If $(\Pi,T,\mySeed)$ does not occur in the 
identified orbit, we reject $\myLatin$ from further consideration.
Second, we test whether $\myLatin$ is the lexicographic minimum of its 
$\Gamma_{\!(\Pi,T,\mySeed)}$-orbit. If not, we reject $\myLatin$ from
further consideration. Otherwise we accept $\myLatin$ as the 
unique representative of its main class. 
See \cite{KO} for a detailed exposition of classification of
combinatorial objects.

From an implementation perspective the present algorithm is
almost identical to the one used in \cite{K,KO2}.
Some additional implementation effort was required to 
work with the groups of order 4, and to implement associated
correctness checks (to be discussed in \S\ref{sect:constructive-cc}). 
In particular, we must list all subgroups of order 4 in 
$\Gamma_{\!\myLatin}$. This was implemented essentially using
brute force, that is, by iterating over all elements of order 4, and 
all pairs of elements of order 2 that generate a group of order 4.
(This strategy suffices because $\Gamma_{\!\myLatin}$ is 
in most cases small; see Table~\ref{tbl:mainclasses}.)

In the search we find 105670178597 Latin squares that extend
the seeds. Among these, we find 1151666641 main classes.
The number of main classes $N_i$ for each $i\geq 3$ is
given later in Table \ref{tbl:mainclasses}.

The search was distributed to a network of 155 Linux PCs using
the batch system \texttt{autoson} \cite{M2}. In total the search 
consumed about 1.3 years of CPU time.

\subsection{Correctness}
\label{sect:constructive-cc}

We carry out a double counting check to gain confidence
in the correctness of the constructive enumeration.
The validation procedures considered in this section and later
in \S\ref{sect:imp} are essential as, for example,
number theoretic properties discussed in \cite{SW}
and elsewhere are not applicable here.

Fix a conjugacy class of admissible groups $[\Pi]$.
We compute in two different ways the total number of tuples
$((\Sigma,T,\mySeed),\myLatin)$ where $(\Sigma,T,\mySeed)$ is 
a seed contained in $\myLatin$ and satisfying $\Sigma\in[\Pi]$.

For a Latin square $\myLatin$, denote by $\seeds_{[\Pi]}(\myLatin)$
the number of seeds associated with $[\Pi]$ contained in $\myLatin$. 
For a seed $(\Sigma,T,\mySeed)$, denote by $\ext(\Sigma,T,\mySeed)$ 
the number of extensions of the seed into a Latin square.
By the Orbit--Stabilizer Theorem, we have
\begin{equation}
\label{eq:constructive-check}
|\Gamma|
\sum_{[\myLatin]}
  \frac{\seeds_{[\Pi]}(\myLatin)}{|\Gamma_{\!\myLatin}|}\,=\ 
|\Gamma|
\sum_{[(\Sigma,T,\mySeed)]:\Sigma\in[\Pi]}
  \frac{\ext(\Sigma,T,\mySeed)}{|\Gamma_{\!(\Sigma,T,\mySeed)}|}\,,
\end{equation}
where the sum on the left-hand side is over the main
classes of Latin squares, and the sum on the right-hand side 
is over the isomorphism classes of seeds associated with $[\Pi]$.

The right-hand side of \eqref{eq:constructive-check} is accumulated
for each classified seed $(\Pi,T,\mySeed)$. In particular, 
$\ext(\Pi,T,\mySeed)$ is simply the number of solutions found in the
exact cover search when extending $(\Pi,T,\mySeed)$.

The left-hand side of \eqref{eq:constructive-check} is accumulated
whenever a constructed $\myLatin$ is accepted as the 
representative of its main class; that is, for every $[\Pi]$ we compute 
$\seeds_{[\Pi]}(\myLatin)$ and accumulate accordingly.
To compute $\seeds_{[\Pi]}(\myLatin)$, we iterate over the
subgroups $\Sigma\leq\Gamma_{\!\myLatin}$ with $\Sigma\in[\Pi]$.
For each such $\Sigma$, we accumulate $\seeds_{[\Pi]}(\myLatin)$
by the value listed in column $M_T$ in Tables~\ref{tbl:primeauto} 
and \ref{tbl:groups4}. To justify this, first observe that whenever
$\Sigma$ and $T$ are fixed, $\myLatin$ determines $\mySeed$ in a seed
$(\Sigma,T,\mySeed)$ occurring in $\myLatin$. In particular, $\mySeed$
is the union of all $\Sigma$-orbits that contain a triple that has
nonempty intersection with $T$. Now, whenever $\Sigma$ is fixed,
$M_T$ counts the number of eligible sets $T$. The values $M_T$ can 
be determined using combinatorial arguments based on the structure
of $\Pi$ and the size and composition constraints for $T$
in relation to $\Pi$. 

Let us give two examples to illustrate 
the combinatorial arguments. First, Type 12 in Table~\ref{tbl:primeauto} 
has $M_T=165$ because there are $3\binom{11}{2}=165$ ways to select 
two distinct points in a common point class. Second, Type 15 in 
Table~\ref{tbl:groups4} has $M_T=9$ because $T$ must consist of 3 
points fixed by $\Sigma$ and must occur as a subset of 
(that is, must be equal to) a triple. 
There are 9 choices for such a $T$ by Lemma~\ref{lem:autotype2}, 
namely the 9 triples of the subsquare of order 3 induced by 
the points fixed by $\Sigma$.

The computed left-hand and right-hand sides of 
\eqref{eq:constructive-check} agree for each $[\Pi]$
in our constructive enumeration, which gives us confidence 
that the results are correct. 
We display the double count values for reference in 
Tables~\ref{tbl:oddprime-check} and \ref{tbl:group4-check}.

\begin{table}
\caption{Double count values for odd prime orders}
\label{tbl:oddprime-check}
\begin{center}
\begin{tabular}{rrrrrrrrr}
\hline
Type &
$p$ & 
$f_R$ & 
$f_C$ & 
$f_S$ & 
$F$ &
Count \\
\hline
10&3   & 2     & 2     & 2     & 3 & 88699187523260511795806208000000      \\
11&3   & 5     & 5     & 5     & 3 & 12309174115893617098752000000000      \\
12&3   & 0     & 0     & 0     & 0 & 19601984696323546934786654208000000   \\
13&5   & 1     & 1     & 1     & 3 & 710224896233056606617600000           \\
14&7   & 4     & 4     & 4     & 3 & 1186258276475008450560000000          \\
15&11  & 0     & 0     & 0     & 3 & 5968708870624483737600000             \\
16&11  & 11    & 0     & 0     & 3 & 38160882055721779200000               \\
\hline
\end{tabular}
\end{center}
\end{table}

\begin{table}
\caption{Double count values for order 4}
\label{tbl:group4-check}
\begin{center}
{\small
\begin{tabular}{rr}
\hline
Type & Count
\\
\hline
2      & 671631524180703313920000000 \\
6      & 2014894572542109941760000000 \\
7      & 3516906890255319171072000000 \\
8      & 54951670160239362048000000000 \\
11     & 1831722338674645401600000000 \\
12     & 5495167016023936204800000000 \\
14     & 21980668064095744819200000000 \\
15     & 307619449557019948744704000000 \\
17     & 39732092772896030588928000000 \\
18     & 13737917540059840512000000000 \\
19     & 268164150381968086794240000000 \\
20     & 65282584150364362113024000000 \\
21     & 179508789190115249356800000000 \\
22     & 16790788104517582848000000000 \\
23     & 1663496959090765967917056000000 \\
25     & 501843074867821516161024000000 \\
26     & 18539743730111373901824000000 \\
28     & 355683843723842644082688000000 \\
29     & 568713151711703904288768000000 \\
30     & 246354443069481975545856000000 \\
31     & 573841974259992911413248000000 \\
32     & 23802621216964125445324800000 \\
33     & 22878993544910822248022016000000 \\
34     & 1477320700587875009298432000000 \\
35     & 464500362123374411907072000000 \\
36     & 12596852214923599074557952000000 \\
37     & 61521447614952423555072000000 \\
38     & 919890958482406920683520000000 \\
39     & 39728022278810086932480000000 \\
42     & 45133638424943262695424000000 \\
43     & 2532319498770923774803968000000 \\
\hline
\end{tabular}
}
\end{center}
\end{table}

\section{Nonconstructive enumeration}
\label{sect:nonconstructive}

\subsection{One-factorizations with symmetry}
\label{sect:onef}

Our objective in this section is to reduce the task of computing 
$|\myddom_{\Pi}|$ for subgroups $\Pi\leq\Gamma$ with $|\Pi|=2$ to 
the task of counting one-factorizations of a complete bipartite
graph with forced symmetry. 

We refer to \cite{W} for basic graph-theoretic terminology.
All graphs considered are undirected, loopless, and without parallel edges.
A {\em one-factor} of a graph $G$ is a spanning $1$-regular subgraph of $G$.
A {\em one-factorization} of $G$ is a set $\mathcal{F}$ of one-factors of $G$ such that
every edge of $G$ occurs in a unique one-factor in $\mathcal{F}$. 
Denote by $\LF(G)$ the number of distinct one-factorizations of $G$.

Let $K_{n,n}$ be the complete bipartite graph with vertex set
$R\cup C$ and bipartition $\{R,C\}$.
Let $\mathcal{F}$ be a one-factorization of $K_{n,n}$, 
and let $f:S\rightarrow\mathcal{F}$ be a bijection that ``labels'' the
one-factors in $\mathcal{F}$ with elements of $S$. 
\begin{Lem}
\label{lem:one-f}
There is a one-to-one correspondence between the 
tuples $(\mathcal{F},f)$ and the Latin squares in $\myddom$. 
\end{Lem}
\begin{Proof}
Each triple $\{r_i,c_j,s_k\}$ 
in a Latin square $\myLatin$ corresponds to the edge $\{r_i,c_j\}$ in 
the one-factor $f(s_k)$ in $\mathcal{F}$. 
\end{Proof}

In particular, Lemma~\ref{lem:one-f} implies
\[
|\myddom|=n!\cdot \LF(K_{n,n})\,.
\]

We now introduce forced symmetry into this setting, that is,
we proceed to study $\myddom_\Pi$ in light of Lemma~\ref{lem:one-f}.
Our task is fortuitously simplified by the assumption 
$|\Pi|=2$, which implies that $\Pi$ fixes at least one of the point classes $R,C,S$
setwise. Without loss of generality (up to conjugation of $\Pi$ in $\Gamma$), 
in what follows we will assume that $\Pi$ fixes $S$ setwise.

Let $\Phi$ be the group consisting of all permutations of $R\cup C$ 
that fix the partition $\{R,C\}$ setwise. Let $\Phi$ act on the set 
of spanning subgraphs of $K_{n,n}$ by permuting the vertices,
and extend the action in the natural way to one-factorizations 
of such graphs.

Observe that each element $\pi\in\Pi$ can be restricted to
$R\cup C$ because $\Pi$ fixes $S$ setwise.
Let $\Delta$ be the restriction of $\Pi$ to $R\cup C$,
and observe that $\Delta\leq\Phi$.

Consider now a Latin square $\myLatin$ with 
$\Pi\leq\Gamma_{\!\myLatin}$, and let $(\mathcal{F},f)$ be the 
corresponding tuple given by Lemma~\ref{lem:one-f}.
We clearly have $\Delta\leq\Phi_{\mathcal{F}}$. However, it is
{\em not} the case that any tuple $(\mathcal{F},f)$ with 
$\Delta\leq\Phi_\mathcal{F}$ corresponds to a Latin square with 
$\Pi\leq\Gamma_{\!\myLatin}$. In particular, the action of $\Delta$ on
$\mathcal{F}$ need not be ``compatible'' with the action of $\Pi$ on $S$.
We proceed with a detailed analysis.

Let $\mathcal{F}$ be a one-factorization of $K_{n,n}$ with 
$\Delta\leq\Phi_\mathcal{F}$, and let $f:S\rightarrow\mathcal{F}$ be a 
bijection. For any $\delta\in\Delta$, define $\bar\delta$ as the permutation 
of $R\cup C\cup \mathcal{F}$ induced by the action of $\delta$ on $R\cup C$ 
and on $\mathcal{F}$. 
(We assume that $R\cup C$ and $\mathcal{F}$ are disjoint sets.)
Note that $\bar\delta$ is well-defined because 
$\delta\in\Delta\leq\Phi_{\mathcal{F}}$. 
Let $\bar\Delta=\{\bar\delta:\delta\in\Delta\}$.
Extend $f$ by the identity mapping on $R\cup C$ to the mapping
$\bar f:R\cup C\cup S\rightarrow R\cup C\cup \mathcal{F}$.
Let 
$\bar f^{-1}\bar\Delta \bar f=
\{\bar f^{-1}\bar\delta \bar f:\delta\in\Delta\}$,
and observe that $\bar f^{-1}\bar\Delta \bar f\leq \Gamma$. 
We say that $(\mathcal{F},f)$ {\em agrees} with $\Pi$ if 
$\bar f^{-1}\bar\Delta \bar f=\Pi$.

\begin{Lem}
\label{lem:agree-pi}
The tuple $(\mathcal{F},f)$ agrees with\/ $\Pi$ if and only if 
the corresponding Latin square $\myLatin$ satisfies 
$\Pi\leq\Gamma_{\!\myLatin}$.
\end{Lem}
\begin{Proof}
Follows from Lemma~\ref{lem:one-f} because $\Pi$ fixes $S$ setwise.
\end{Proof}

We say that $\mathcal{F}$ is {\em compatible} with $\Pi$ if there exists 
an $f:S\rightarrow\mathcal{F}$ such that $(\mathcal{F},f)$ agrees with $\Pi$.
Denote by $\Psi$ the group of all permutations of $R\cup C\cup S$ that fix
$R\cup C$ pointwise. 
\begin{Lem}
\label{lem:num-agree}
Let $\mathcal{F}$ be compatible with\/ $\Pi$. Then there are
exactly $|\Psi_\Pi|$ bijections $f$ such that $(\mathcal{F},f)$ 
agrees with $\Pi$.
\end{Lem}
\begin{Proof}
Fix a reference bijection and establish a one-to-one
correspondence between the set of all bijections and $\Psi_\Pi$.
\end{Proof}

Denote by $\LF(K_{n,n},\Pi)$ the number of distinct 
one-factorizations of $K_{n,n}$ compatible with $\Pi$. 
From Lemmata~\ref{lem:agree-pi} and \ref{lem:num-agree} it thus follows
that
\begin{equation}
\label{eq:l-pi-size}
|\myddom_\Pi|=|\Psi_\Pi|\cdot \LF(K_{n,n},\Pi)\,.
\end{equation}

The next three lemmata rely on our assumption that $\Pi$ has order 2.

\begin{Lem}
Suppose $\Pi$ restricted to $S$ has 
$a_1$ orbits of length\/ $1$ and $a_2$ orbits of length\/ $2$. 
Then 
\begin{equation}
\label{eq:psi-pi}
|\Psi_\Pi|=a_1!a_2!2^{a_2}\,.
\end{equation}
\end{Lem}
\begin{Proof}
There are $a_1!a_2!2^{a_2}$ permutations in $\Psi$
that fix the nonidentity element of $\Pi$ under conjugation. 
\end{Proof}

\begin{Lem}
\label{lem:compat-one-f}
Let $\mathcal{F}$ be a one-factorization of $K_{n,n}$ with 
$\Delta\leq\Phi_\mathcal{F}$.
Then $\mathcal{F}$ is compatible with\/ $\Pi$ if and only if
the action of\/ $\Delta$ on $\mathcal{F}$ 
has an equal number of orbits of each length
as the action of\/ $\Pi$ on $S$.
\end{Lem}
\begin{Proof}
The permutation of $S$ induced by the nonidentity element of $\Pi$ 
can be relabelled (by conjugation with an appropriate $f$) to 
the permutation of $\mathcal{F}$ induced by the nonidentity element 
of $\Delta$ if and only if there are an equal number of 
cycles of each length.
\end{Proof}

Let $G$ be a spanning subgraph of $K_{n,n}$.
We say that a one-factorization $\mathcal{F}$ of 
$G$ is a $(\Delta,\vec a)$-{\em factorization} if 
(a) $\Delta\leq\Phi_{\mathcal{F}}$; and 
(b) for each $i=1,2$, the action of $\Delta$ on $\mathcal{F}$ has 
$a_i$ orbits of length $i$ with $\vec a=(a_1,a_2)$.
Denote by $\LFp{\Delta}{\vec a}{G}$ the number of distinct
$(\Delta,\vec a)$-factorizations of $G$.

\begin{Lem}
Suppose $\Pi$ restricted to $S$ has 
$a_1$ orbits of length\/ $1$ and $a_2$ orbits of length\/ $2$. 
Then 
\begin{equation}
\label{eq:delta-pi}
\LF(K_{n,n},\Pi)=\LFp{\Delta}{\vec a}{K_{n,n}}\,.
\end{equation}
\end{Lem}
\begin{Proof}
Immediate from Lemma~\ref{lem:compat-one-f}.
\end{Proof}

We have now reduced the task of computing
$|\myddom_\Pi|$ via 
\eqref{eq:l-pi-size}, \eqref{eq:psi-pi}, and \eqref{eq:delta-pi}
to computing $\LFp{\Delta}{\vec a}{K_{n,n}}$.

We employ two different methods to
compute the values $\LFp{\Delta}{\vec a}{K_{n,n}}$ for
the conjugacy classes of order 2 in Table~\ref{tbl:primeauto}.
The first method stems from a recursion over regular
spanning subgraphs of $K_{n,n}$ that has evolved over 
a period of more than half a century into the form 
presented by McKay and Wanless \cite{MW}. 
The second method extends a ``forward accumulation'' technique 
used 
in \cite{KO2} to 
count the number of distinct one-factorizations of $K_{14}$.

Our present contribution is to modify
both methods to take into account the symmetry
forced by $\Delta$ and the required structure $\vec a=(a_1,a_2)$ 
for the lengths of the $\Delta$-orbits of one-factors.

\subsection{Backward recursion}
\label{sect:bwd-rec}

For convenience in what follows, let us 
abbreviate $\vec e_1=(1,0)$ and $\vec e_2=(0,1)$.
We say that a spanning subgraph $G$ of $K_{n,n}$ is 
$(\Delta,\vec a)$-{\em factorizable} if $G$ has at least 
one $(\Delta,\vec a)$-factorization.
Observe that a $(\Delta,\vec a)$-factorizable 
graph is necessarily $k$-regular with $k=\sum_i ia_i$. 
Furthermore, each $(\Delta,\vec a)$-factorizable graph decomposes 
into an edge-disjoint union consisting of $a_1$ $(\Delta,\vec e_1)$-factorizable
and $a_2$ $(\Delta,\vec e_2)$-factorizable graphs.
Denote by $\onef G\Delta i$ the set of all
$(\Delta,\vec e_i)$-factorizable spanning subgraphs of $G$.

Counting in two different ways the number of distinct 
$(\Delta,\vec a)$-factorizations of $G$ with
one individualized $\Delta$-orbit consisting 
of $i$ one-factors, we obtain, for $a_i\geq 1$,
\begin{equation}
\label{eq:lf}
a_i\cdot\LFp{\Delta}{\vec a}{G}=
\sum_{F\in\onef G \Delta i}\LFp{\Delta}{\vec e_i}{F}\cdot\LFp{\Delta}{\vec a-\vec e_i}{G-F}\,.
\end{equation}

In particular, if we know (for example, by recursion) the value
$\LFp{\Delta}{\vec a-\vec e_i}{H}$ for all 
$(\Delta,\vec a-\vec e_i)$-factorizable $H$, 
we can use \eqref{eq:lf} to compute 
$\LFp{\Delta}{\vec a}{G}$ for all 
$(\Delta,\vec a)$-factorizable $G$.
Furthermore, $\LFp{\Delta}{\vec a}{G}$ needs to be computed
for only one representative $(G,\Delta)$ chosen from the 
$\Phi$-orbit $[(G,\Delta)]$. Here $\Phi$ acts on $G$ by 
permutation of vertices and on $\Delta$ by conjugation.

This is analogous to the method (without the forced $\Delta$ and $\vec a$) 
used in \cite{MMM} and other earlier studies on counting Latin squares.

\subsection{Forward accumulation}
\label{sect:fwd-acc}

We can visit every $(\Delta,\vec a)$-factorizable orbit $[(G,\Delta)]$
with $a_i\geq 1$ by the following procedure:
for each $(\Delta,\vec a-\vec e_i)$-factorizable orbit $[(H,\Delta)]$,
consider exactly one tuple $(H,\Delta)$ from the orbit;
for each such tuple $(H,\Delta)$, consider each graph
$F\in\onef {K_{n,n}-H}\Delta i$; for each tuple $(H,F,\Delta)$,
visit the $(\Delta,\vec a)$-factorizable orbit $[(H\cup F,\Delta)]$.

To compute the value $\LFp{\Delta}{\vec a}{G}$ for each visited 
orbit $[(G,\Delta)]$, we associate with $[(G,\Delta)]$ an accumulator
variable $x_{[(G,\Delta)]}$ that is initially set to zero
and incremented whenever $[(G,\Delta)]$ is visited. 
Our objective is to have the value $a_i\cdot\LFp{\Delta}{\vec a}{G}$ 
in $x_{[(G,\Delta)]}$ when the procedure terminates.

We proceed to analyze the visiting procedure, with the objective
of determining an appropriate increment to the counter variable
on each visit. The following two lemmata are immediate consequences
of the Orbit-Stabilizer Theorem (cf.\ \mbox{\cite[Lemmata 1,2]{KO2}).}

\begin{Lem}
\label{lem:decomp}
Any tuple $(G,\Delta)$ in the orbit\/ $[(H\cup F,\Delta)]$ admits exactly 
\[
\sigma(H,F,\Delta)=\frac{|\Phi_{(H\cup F,\Delta)}|}{|\Phi_{(H,F,\Delta)}|}
\]
decompositions $G=H'\cup F'$ into tuples $(H',F',\Delta)$ 
in the orbit\/ $[(H,F,\Delta)]$.
\end{Lem}

\begin{Lem}
\label{lem:visits}
The procedure visits an orbit\/ $[(H\cup F,\Delta)]$ exactly 
\[
\tau(H,F,\Delta)=\frac{|\Phi_{(H,\Delta)}|}{|\Phi_{(H,F,\Delta)}|}
\]
times via tuples $(H',F',\Delta)$ in the orbit\/ $[(H,F,\Delta)]$.
\end{Lem}

It now follows from Lemma \ref{lem:decomp} 
and \eqref{eq:lf} that, for $a_i\geq 1$, 
\begin{equation}
\label{eq:dgm-up}
\begin{split}
a_i\cdot \LFp{\Delta}{\vec a}{G}
&=\sum_{(H,F,\Delta):H\cup F=G}\LFp{\Delta}{\vec e_i}{F}\cdot\LFp{\Delta}{\vec a-\vec e_i}{H}\\
&=\sum_{[(H,F,\Delta)]:H\cup F=G}\sigma(H,F,\Delta)\cdot\LFp{\Delta}{\vec e_i}{F}\cdot\LFp{\Delta}{\vec a-\vec e_i}{H}\,.
\end{split}
\end{equation}

This observation enables us to accumulate the value 
$a_i\cdot\LFp{\Delta}{\vec a}{G}$ to $x_{[(G,\Delta)]}$.
Namely, each time $[(G,\Delta)]$ 
is visited via a tuple $(H,F,\Delta)$, 
we increment $x_{[G,\Delta)]}$ by the rule
\begin{equation}
\label{eq:acc}
x_{[(G,\Delta)]}\leftarrow x_{[(G,\Delta)]}+
                        \frac{\sigma(H,F,\Delta)}{\tau(H,F,\Delta)}\cdot\LFp{\Delta}{\vec e_i}{F}\cdot\LFp{\Delta}{\vec a-\vec e_i}{H}\,.
\end{equation}
Equivalently, for each tuple $(H,F,\Delta)$ considered by the
visiting procedure, we apply the rule
\begin{equation}
\label{eq:acc2}
x_{[(H\cup F,\Delta)]}\leftarrow x_{[(H\cup F,\Delta)]}+
                        \frac{|\Phi_{(H\cup F,\Delta)}|}{|\Phi_{(H,\Delta)}|}\cdot\LFp{\Delta}{\vec e_i}{F}\cdot\LFp{\Delta}{\vec a-\vec e_i}{H}\,.
\end{equation}

Analogously to \cite[Lemma 3]{KO2}, we have the
following result.

\begin{Lem}
The total accumulation to $x_{[(G,\Delta)]}$ is $a_i\cdot \LFp{\Delta}{\vec a}{G}$.
\end{Lem}
\begin{Proof}
By Lemma \ref{lem:visits} and \eqref{eq:acc}, 
the total accumulation to $x_{[(G,\Delta)]}$ from an orbit $[(H,F,\Delta)]$ satisfying 
$G=H\cup F$ is $\sigma(H,F,\Delta)\cdot\LFp{\Delta}{\vec e_i}{F}\cdot\LFp{\Delta}{\vec a-\vec e_i}{H}$. 
Taking the sum over all such classes, the claim follows by \eqref{eq:dgm-up}.
\end{Proof}

\subsection{Implementation and correctness}

\label{sect:imp}
Associated with each graph $F$ in \S\ref{sect:bwd-rec} and \S\ref{sect:fwd-acc}
there are exactly $\LFp{\Delta}{\vec e_i}{F}$ distinct $\Delta$-orbits 
(of length $i$) of one-factors; that is, the distinct 
$(\Delta,\vec e_i)$-factorizations of $F$. From an implementation
perspective it is convenient to iterate over these $\Delta$-orbits 
directly instead of the underlying graphs $F$. In particular, we achieve the 
required accumulation in \eqref{eq:lf} and \eqref{eq:acc2} without 
computing the values $\LFp{\Delta}{\vec e_i}{F}$ explicitly.

Starting with the empty spanning subgraph of $K_{n,n}$, 
we use the forward accumulation method to compute the values
$\LFp{\Delta}{\vec a}{G}$ for all $(\Delta,\vec a)$-factorizable
orbits $[(G,\Delta)]$. 
We carry out the forward accumulation in two different ways. 
First, we accumulate the $\Delta$-orbits in order of increasing length; 
that is, we first add all $\Delta$-orbits of length $1$, and then 
the $\Delta$-orbits of length $2$. Second, as a consistency check, 
we add the orbits in order of decreasing length.

In addition to the forward accumulation method, 
we verify the results using the recursion \eqref{eq:lf}, that is,
when we have accumulated the value $\LFp{\Delta}{\vec a}{G}$ by 
addition of orbits of length $i$, we apply \eqref{eq:lf} to check the
result using the values $\LFp{\Delta}{\vec a-\vec e_i}{H}$.

The two different orders of adding orbits enable yet another 
way to cross-check the results (cf.~\cite{DGM} and \cite{KO2}).
Namely, 

\begin{equation}
\label{eq:mim}
\LFp{\Delta}{\vec a+\vec b}{K_{n,n}}=
\frac{%
\sum_{[(G,\Delta)]}\frac{|\Phi_{\Delta}|}{|\Phi_{(G,\Delta)}|}
\cdot\LFp{\Delta}{\vec a}{G}\cdot\LFp{\Delta}{\vec b}{K_{n,n}-G}
}{\binom{a_1+b_1}{a_1}\binom{a_2+b_2}{a_2}}
\end{equation}

\noindent
holds for every nonnegative $\vec a=(a_1,a_2)$ 
and $\vec b=(b_1,b_2)$ with $\sum_i i(a_i+b_i)=n$,
where the sum is over all $(\Delta,\vec a)$-factorizable
orbits $[(G,\Delta)]$.
Due to the two different orders of adding orbits, we cannot
utilise all possible decompositions $\vec a+\vec b$ for checking, 
however. For example, if we add orbits in the orders 1,1,1,2,2,2,2
and 2,2,2,2,1,1,1, then we can use eight different decompositions
in \eqref{eq:mim}.

The approaches were implemented using the C programming language
with the help of three software libraries:
\emph{nauty} \cite{M1} for isomorphism testing,
the GNU Multiple Precision Arithmetic Library \cite{gmplib}
for handling large integers (and intermediate rationals),
and \emph{libexact} \cite{KP} for simplifying the 
search for $\Delta$-orbits of one-factors. 
The number of orbits $[(G,\Delta)]$ at any step of the search 
never exceeded 200000, so memory requirement was 
not a major issue of the search (compare with \cite{KO2,MMM}).
The computations took about 20 days on a contemporary Linux PC.

\subsection{Results}

The numbers of $(\Delta,\vec a)$-factorizable orbits $[(G,\Delta)]$
are shown in Tables \ref{tab:int1} and \ref{tab:int2} for
situations when the one-factor orbits are appended in order
of increasing and decreasing length $i$, respectively.
The leftmost column contains the number of one-factor orbits $\sum_i a_i$,
and there is one column for each of the conjugacy
classes of order $2$ in Table~\ref{tbl:primeauto}.

\begin{table}[htbp]
\caption{Number of $(\Delta,\vec a)$-factorizable orbits, increasing
length}
\label{tab:int1}
\begin{center}
\begin{tabular}{rrrrrrrrrr}\hline
\raisebox{-0.5mm}{$\sum_i a_i$ $\backslash$ Type\hspace*{-1.5cm}}&
\raisebox{-0.5mm}{\hspace*{1.5cm}$1$}&\raisebox{-0.5mm}{$2$}&
\raisebox{-0.5mm}{$3$}&\raisebox{-0.5mm}{$4$}&
\raisebox{-0.5mm}{$5$}&\raisebox{-0.5mm}{$6$}&
\raisebox{-0.5mm}{$7$}&\raisebox{-0.5mm}{$8$}&
\raisebox{-0.5mm}{$9$}\\[0.5mm]\hline
0 & 1   &1   &1  &1     &1   &1     &1     &1     &1     \\
1 & 1  	&1   &1  &6     &6   &6     &6     &6     &6     \\ 
2 & 304	&10  &10 &1089  &22  &22    &22    &22    &22    \\ 
3 &34792&16  &10 &170321&1365&1365  &1365  &1365  &1365  \\ 
4 &8530	&949 &5  &33851 &184239&30429&30429 &30429 &30429 \\ 
5 &20  	&1464&1  &46    &35286&146085&146085&146085&146085\\ 
6 & 1  	&11  &1  &1     &49  &32039 &125763&125763&125763\\ 
7 &    	&1   &1  &      &1   &45    &19567 &19567 &19567 \\ 
8 &    	&    &1  &      &    &1     &36    &630   &630   \\ 
9 &     &    &   &      &    &      &1     &21    &21    \\ 
10&     &    &   &      &    &      &      &1     &1     \\ 
11&     &    &   &      &    &      &      &      &1     \\\hline
\end{tabular}
\end{center}
\end{table}

\begin{table}[htbp]
\caption{Number of $(\Delta,\vec a)$-factorizable orbits, decreasing length}
\label{tab:int2}
\begin{center}
\begin{tabular}{rrrrrrrrrr}\hline
\raisebox{-0.5mm}{$\sum_i a_i$ $\backslash$ Type\hspace*{-1.5cm}}&
\raisebox{-0.5mm}{\hspace*{1.5cm}$1$}&\raisebox{-0.5mm}{$2$}&
\raisebox{-0.5mm}{$3$}&\raisebox{-0.5mm}{$4$}&
\raisebox{-0.5mm}{$5$}&\raisebox{-0.5mm}{$6$}&
\raisebox{-0.5mm}{$7$}&\raisebox{-0.5mm}{$8$}&
\raisebox{-0.5mm}{$9$}\\[0.5mm]\hline
0 & 1   &1   &1  &1     &1   &1     &1     &1     &1     \\
1 & 5  	&6   &1  &6     &6   &6     &6     &6     &6     \\ 
2 &8528	&1464&3  &266   &266 &266   &266   &1089  &22    \\ 
3 &34792&1182&1  &266   &266 &266   &170321&35186 &1365  \\ 
4 &319 	&16  &5  &6     &6   &33851 &187719&184239&30429 \\ 
5 & 1  	&10  &10 &1     &46  &1524  &35286 &178732&146085\\ 
6 & 1  	&1   &10 &1     &5   &49    &1475  &32039 &125763\\ 
7 &    	&1   &1  &      &1   &4     &45    &1220  &19567 \\ 
8 &    	&    &1  &      &    &1     &3     &36    &630   \\ 
9 &     &    &   &      &    &      &1     &2     &21    \\ 
10&     &    &   &      &    &      &      &1     &1     \\ 
11&     &    &   &      &    &      &      &      &1     \\\hline
\end{tabular}
\end{center}
\end{table}

The second column of Table \ref{tbl:res2} summarizes the results of 
the search, which enable us to determine the values $|\myddom_{[\Pi]}|$ 
in the rightmost column of Table \ref{tbl:res2} via the combination 
of \eqref{eq:lpi-size}, \eqref{eq:l-pi-size},
\eqref{eq:psi-pi}, and \eqref{eq:delta-pi}. That is,
\[
|\myddom_{[\Pi]}| =
\frac{|\Gamma|}{|\Gamma_{\!\Pi}|}
\cdot a_1!a_2!2^{a_2}
\cdot \LFp{\Delta}{(a_1,a_2)}{K_{11,11}},
\]
where $|\Gamma|=3!(11!)^3$ and the values 
$|\Gamma_{\!\Pi}|$ are displayed in Table~\ref{tbl:gammapi}; 
these are derived by considering the number of permutations in $\Gamma$ 
that fix the nonidentity element of $\Pi$ under conjugation.

\begin{table}
\caption{$|\Gamma_{\!\Pi}|$ for conjugacy classes with $|\Pi|=2$}
\label{tbl:gammapi}
\begin{center}
\begin{tabular}{rr@{\,}r}\hline
Type & & $|\Gamma_{\!\Pi}|$ \\
\hline
1 & $(1!\cdot 5!\cdot 2^{5})^3\cdot 3!=$   &      339738624000\\
2 & $(3!\cdot 4!\cdot 2^{4})^3\cdot 3!=$   &       73383542784\\
3 & $(5!\cdot 3!\cdot 2^{3})^3\cdot 3!=$   &     1146617856000\\
4 & $1!\cdot 5!\cdot 2^5\cdot 11!\cdot 2=$ &      306561024000\\
5 & $3!\cdot 4!\cdot 2^4\cdot 11!\cdot 2=$ &      183936614400\\
6 & $5!\cdot 3!\cdot 2^3\cdot 11!\cdot 2=$ &      459841536000\\
7 & $7!\cdot 2!\cdot 2^2\cdot 11!\cdot 2=$ &     3218890752000\\
8 & $9!\cdot 1!\cdot 2^1\cdot 11!\cdot 2=$ &    57940033536000\\
9 & $\quad11!\cdot 0!\cdot 2^0\cdot 11!\cdot 2=$&  3186701844480000\\
\hline
\end{tabular}
\end{center}
\end{table}

\begin{table}
\caption{$|\myddom_{[\Pi]}|$ for conjugacy classes with $|\Pi|=2$}
\label{tbl:res2}
\begin{center}
\begin{tabular}{rrr}\hline
Type &  $\LF(K_{11,11},\Delta,\vec a)$ & $|\myddom_{[\Pi]}|$\\\hline
1 &   157811617463135109120& 680681075601465561779007809126400000\\
2 &    16482860057738870784& 197485300889612684060962848768000000\\
3 &       57076088832000000&    109414881030165485322240000000000\\
4 &       81867724734136320&    391332044120479692845639270400000\\
5 &    91614680894074060800& 437923008879888637579625496576000000\\
6 &   984580848455978188800&4706348408721523581497414516736000000\\
7 &   588014650532826316800&2810741057151330357801067216896000000\\
8 &    39448015149028147200& 188563593954724468142515421184000000\\
9 &      252282619805368320&   1205924234796020705924638310400000\\\hline
\end{tabular}
\end{center}
\end{table}

\section{Number of equivalence classes}
\label{sect:completing}

\subsection{Main classes}
As a result of the computations described in \S\ref{sect:constructive}
and \S\ref{sect:nonconstructive}, we have a constructive enumeration
for all main classes with an autoparatopy group of order at least 3.
The numbers of such main classes are given in Table~\ref{tbl:mainclasses}
for $i\geq 3$. Table \ref{tbl:res2} enables us to solve the number of
main classes with an autoparatopy group of order 2. Indeed, we 
compute $|\mydddom_{[\Pi]}|$ by direct summation over 
the constructed main classes \eqref{eq:lc}, and use \eqref{eq:main}
together with Table~\ref{tbl:res2} to arrive at Table~\ref{tbl:ord2}.

\begin{table}
\caption{$|\mydddom_{[\Pi]}|$ and $N_{[\Pi]}$ 
         for conjugacy classes with $|\Pi|=2$}
\label{tbl:ord2}
\begin{center}
\begin{tabular}{rrr}
\hline
Type & $|\mydddom_{[\Pi]}|$ & $N_{[\Pi]}$\\
\hline
1& 1788612130139530711950950400000 & 3567419044431 \\
2& 2090765465832065147338752000000 & 1035003382971 \\
3& 40157841013697699905536000000   & 573229534 \\
4& 37761778205327560959590400000   & 2050761205 \\
5& 1253232559784675856089088000000 & 2295134347303 \\
6& 3238510153755989660663808000000 & 24665809156818 \\
7& 2028400471919270993854464000000 & 14730996126592 \\
8& 252164628166875895037952000000  & 988254629517 \\
9& 14823428114332518344294400000   & 6320133846 \\
\hline
$N_2$ &                              & 47291560812217\\
\hline
\end{tabular}
\end{center}
\end{table}

Table \ref{tbl:ord2} and Table \ref{tbl:mainclasses} (for $i\geq 3$) 
give
\begin{equation}
\label{eq:totaln2}
|\Gamma|\sum_{i\geq 2} \frac{N_i}{i}=
9023583561995938862980803959193600000\,.
\end{equation}

The number of distinct Latin squares of order 11 is \cite{MW}
\begin{equation}
\label{eq:all-label}
\begin{split}
|\mydom|&=
11!\cdot 10!\cdot 5363937773277371298119673540771840\\
&=776966836171770144107444346734230682311065600000\,.
\end{split}
\end{equation}
Solving for $N_1$ in 
\eqref{eq:totaln} using \eqref{eq:totaln2} and \eqref{eq:all-label}, we obtain
\[
N_1=2036029552535590421717241
\]
and hence
\[
\sum_i N_i=2036029552582883134196099\,.
\]

Table~\ref{tbl:mainclasses} gives the number $N_i$ of main classes
for each order $i$ of the autoparatopy group.

\begin{table}
\caption{Main classes of Latin squares of order $11$}
\label{tbl:mainclasses}
\begin{center}
\begin{tabular}{rr}
\hline
$i$ & $N_i$\\
\hline
1 &2036029552535590421717241\\
2 &47291560812217\\
3 &1111651266\\
4 &39004721\\
5 &9131\\
6 &960771\\
7 &1294\\
8 &30390\\
9 &2636\\
10 &323\\
11 &3\\
12 &4191\\
14 &105\\
15 &4\\
16 &631\\
18 &625\\
20 &37\\
21 &37\\
22 &4\\
\hline
\end{tabular}
\ \ \ 
\begin{tabular}{rr}
\hline
$i$ & $N_i$\\
\hline
24 &274\\
27 &27\\
30 &4\\
32 &16\\
36 &86\\
40 &4\\
42 &5\\
48 &9\\
54 &18\\
60 &1\\
63 &3\\
72 &14\\
108 &5\\
110 &2\\
120 &1\\
162 &1\\
324 &1\\
7260 &1\\
\hline
$\sum_i N_i$ &2036029552582883134196099\\
\hline
\end{tabular}
\end{center}
\end{table}

\subsection{One-factorizations and isotopy classes}

For a Latin square $\myLatin$, denote by $\Par(\myLatin)$ the 
autoparatopy group of $\myLatin$, and denote by $\Is(\myLatin)$ the
subgroup of $\Par(\myLatin)$ that fixes each of the point classes $R,C,S$ setwise.
This group is called the {\em autotopy} group of $\myLatin$, and 
its elements are called {\em autotopisms}. 
Let $\Ty(\myLatin)=3,2,1,1$ when 
$|\Par(\myLatin)|/|\Is(\myLatin)|=1,2,3,6$, respectively.
Observe that these quantities are independent of the main class
representative $\myLatin$ because $\Is(\myLatin)$ is a normal subgroup 
of $\Par(\myLatin)$.

By \cite[Theorem 2(ii) and 2(iii)]{MMM}, each main class
$[\myLatin]$ splits into $\Ty(\myLatin)$ isomorphism classes 
of one-factorizations of $K_{n,n}$ (cf.~\S\ref{sect:onef})
and into 
\mbox{$3!\cdot|\Is(\myLatin)|/|\Par(\myLatin)|$} isotopy classes.

Taking the sum over all main classes, we obtain 
Theorem~\ref{thm:main}(ii) and \ref{thm:main}(iii).
Table~\ref{tbl:paris} gives the order of the autoparatopy group and
the order of the autotopy group for each main class.

\begin{table}
\caption{Orders of autoparatopy and autotopy groups}
\label{tbl:paris}
\begin{tabular}{rrr}
\hline
$|\Par|$ & $|\Is|$ & Main classes\\
\hline
1 & 1 & 2036029552535590421717241\\
2 & 1 & 42688565155281\\
2 & 2 & 4602995656936\\
3 & 1 & 933551378\\
3 & 3 & 178099888\\
4 & 2 & 36528967\\
4 & 4 & 2475754\\
5 & 5 & 9131\\
6 & 1 & 742503\\
6 & 2 & 28147\\
6 & 3 & 164658\\
6 & 6 & 25463\\
7 & 7 & 1294\\
8 & 4 & 27705\\
8 & 8 & 2685\\
9 & 3 & 2125\\
9 & 9 & 511\\
10 & 5 & 275\\
10 & 10 & 48\\
11 & 11 & 3\\
12 & 2 & 2201\\
12 & 4 & 400\\
12 & 6 & 1470\\
12 & 12 & 120\\
14 & 7 & 79\\
14 & 14 & 26\\
15 & 5 & 4\\
16 & 8 & 607\\
16 & 16 & 24\\
18 & 3 & 486\\
18 & 6 & 35\\
18 & 9 & 96\\
\hline
\end{tabular}
\ \ \ \ 
\begin{tabular}{rrr}
\hline
$|\Par|$ & $|\Is|$ & Main classes\\
\hline
18 & 18 & 8\\
20 & 10 & 22\\
20 & 20 & 15\\
21 & 7 & 32\\
21 & 21 & 5\\
22 & 11 & 2\\
22 & 22 & 2\\
24 & 4 & 213\\
24 & 8 & 6\\
24 & 12 & 55\\
27 & 9 & 27\\
30 & 5 & 2\\
30 & 10 & 2\\
32 & 16 & 16\\
36 & 6 & 70\\
36 & 12 & 10\\
36 & 18 & 6\\
40 & 20 & 4\\
42 & 14 & 2\\
42 & 21 & 3\\
48 & 8 & 9\\
54 & 9 & 17\\
54 & 27 & 1\\
60 & 10 & 1\\
63 & 21 & 3\\
72 & 12 & 14\\
108 & 18 & 5\\
110 & 55 & 2\\
120 & 60 & 1\\
162 & 27 & 1\\
324 & 54 & 1\\
7260 & 1210 & 1\\
\hline
\end{tabular}
\end{table}

\subsection{Quasigroups and loops}

A permutation with $a_i$ cycles of length $i$ for
$i \geq 1$ is said to have cycle structure $(a_1,a_2,\ldots)$.
For an autotopism $\alpha\in\Is(\myLatin)$, define
$f(\alpha)=\prod_i a_i!\,i^{a_i}$ and $g(\alpha)=a_1$
if $\alpha$ is not the identity permutation and has the 
same cycle structure $(a_1,a_2,\ldots)$ in each of the point classes $R,C,S$; 
otherwise $f(\alpha)=g(\alpha)=0$.

The number of isomorphism classes of quasigroups is,
by \cite[Theorem~4]{MMM},
\begin{equation}
\label{eq:qg}
\frac{|\mydom|}{n!}+
\sum_{[\myLatin]}\frac{3!}{|\Par(\myLatin)|}\sum_{\alpha\in\Is(\myLatin)}f(\alpha)^2\,.
\end{equation}
The number of isomorphism classes of loops is,
by \cite[Theorem~5]{MMM},
\begin{equation}
\label{eq:loop}
\frac{|\mydom|}{n!\cdot(n-1)!\,^2}+
\sum_{[\myLatin]}\frac{3!}{|\Par(\myLatin)|}\sum_{\alpha\in\Is(\myLatin)}g(\alpha)^2\,.
\end{equation}

We obtain Theorem~\ref{thm:main}(iv) and \ref{thm:main}(v)
by \eqref{eq:all-label}, \eqref{eq:qg}, \eqref{eq:loop}, 
and summation over the main classes.

\section*{Acknowledgments}

The authors thank Ian Wanless and an
anonymous referee for valuable comments and corrections.
This research was supported in part by 
the Academy of Finland, Grants 107493, 110196, 117499, 130142,
and 132122; 
as well as by the National Science Foundation, Grant Number 0633333.


\end{document}